\documentclass{amsart}

\usepackage{amssymb}

\allowdisplaybreaks

\newtheorem{theorem}{Theorem}[section]
\newtheorem{lemma}[theorem]{Lemma}

\newtheorem{corollary}[theorem]{Corollary}
\newtheorem{remark}[theorem]{Remark}

\newtheorem{maintheo}{Theorem}
\newtheorem{theoA}{Theorem}
\newtheorem{theoB}{Theorem}
\newtheorem{theoC}{Theorem}

\theoremstyle{definition}

\newcommand{\R}{\mathbb{R}}
\newcommand{\C}{\mathbb{C}}

\newcommand{\la}{\mathbb{\lambda}}

\newcommand{\ten}{\otimes}

\newcommand{\dem}{\noindent {\bf Proof. }}

\newcommand{\demC}{\noindent {\bf Proof of Theorem \ref{QSQSQS}. }}
\newcommand{\fin}{\hspace*{\fill} $\square$ \vskip0.2cm}
\newcommand{\prodd}{\prod\nolimits}

\def\bubl{{\displaystyle \mathop{\mathsf{A}}^{\circ}}}

\newcommand{\summ}{\sum\nolimits}

\begin{document}

\title[Operator space $L_p$ embedding theory I]
{Operator space $L_p$ embedding theory I}

\author[Junge and Parcet]
{Marius Junge$^{\ast}$ and Javier Parcet$^{\dag}$}

\footnote{$^{\ast}$Partially supported by the NSF DMS-0301116.}
\footnote{$^{\dag}$Partially supported by \lq Programa Ram{\'o}n y
Cajal, 2005\rq${}$ and \\ \indent also by Grants MTM2004-00678 and
CCG06-UAM/ESP-0286, Spain.} 

\begin{abstract}
Let $\mathrm{X}_1$ and $\mathrm{X}_2$ be subspaces of quotients of
$R \oplus \mathrm{OH}$ and $C \oplus \mathrm{OH}$ respectively. We
use new free probability techniques to construct a completely
isomorphic embedding of the Haagerup tensor product $\mathrm{X}_1
\otimes_h \mathrm{X}_2$ into the predual of a sufficiently large
QWEP von Neumann algebra. As an immediate application, given any
$1 < q \le 2$, our result produces a completely isomorphic
embedding of $\ell_q$ (equipped with its natural operator space
structure) into $L_1(\mathcal{A})$ with $\mathcal{A}$ a QWEP von
Neumann algebra.
\end{abstract}

\maketitle

\vskip-1cm

\null




\section*{Introduction}

The idea of replacing functions by linear operators, the process
of quantization, goes back to the foundations of quantum mechanics
and has a great impact in mathematics. This applies for instance
to representation theory, noncommutative geometry, operator
algebra, quantum and free probability or operator space theory
among other fields. The quantization of measure theory leads to
the theory of $L_p$ spaces defined over general von Neumann
algebras, so called \emph{noncommutative $L_p$ spaces}. This
theory was initiated by Segal, Dixmier and Kunze in the fifties
and continued years later by Haagerup, Fack, Kosaki and many
others. We refer to the recent survey \cite{PX2} for a complete
exposition. In a series of papers beginning with this work, we
will investigate noncommutative $L_p$ spaces in the language of
noncommutative Banach spaces, so called \emph{operator spaces}.

\vskip5pt

The theory of operator spaces took off in 1988 with Ruan's work
\cite{R} and it has been developed since then by Blecher/Paulsen,
Effros/Ruan and Pisier as a noncommutative generalization of
Banach space theory \cite{ER,Pau,P3}. In his book \cite{P2} on
vector-valued noncommutative $L_p$ spaces, Pisier considered a
distinguished operator space structure on $L_p$. In fact, the
right category when dealing with noncommutative $L_p$ is in many
aspects that of operator spaces. Indeed, this has become clear in
the last years by recent results on noncommutative martingales and
related topics. In this and forthcoming papers we shall prove a
fundamental structure theorem of $L_p$ spaces in the category of
operator spaces, solving a problem formulated by Gilles Pisier.

\begin{maintheo} Let $1 \le p < q \le 2$ and let $\mathcal{M}$ be a
von Neumann algebra. Then, there exists a sufficiently large von
Neumann algebra $\mathcal{A}$ and a completely isomorphic
embedding of $L_q(\mathcal{M})$ into $L_p(\mathcal{A})$, where
both spaces are equipped with their respective natural operator
space structures. Moreover, we have
\begin{itemize}
\item[\textbf{i)}] If $\mathcal{M}$ is $\mathrm{QWEP}$, we can
choose $\mathcal{A}$ to be $\mathrm{QWEP}$.

\item[\textbf{ii)}] If $\mathcal{M}$ is hyperfinite, we can choose
$\mathcal{A}$ to be hyperfinite.

\item[\textbf{iii)}] If $\mathcal{M}$ is infinite-dimensional,
then $\mathcal{A}$ must be of type $\mathrm{III}$.
\end{itemize}
\end{maintheo}

In the category of operator spaces, the $L_p$ embedding theory
presents some significant differences, some of them already apply to
$\mathrm{OH}$. For example,
in sharp contrast with the classical theory, it was proved in
\cite{J00} that $\mathrm{OH}$ does not embed completely
isomorphically into any $L_p$ space for $2 < p < \infty$.
Moreover, after \cite{P5} we know that there is no possible
cb-embedding of $\mathrm{OH}$ into the predual of a semifinite von
Neumann algebra. As it follows from our result, this also happens
for $1 \le p < q \le 2$ and justifies the relevance of type III
von Neumann algebras in the subject. In this paper we prove only a
part of the Theorem stated above. Our concern in this paper will
be to prove a particular case of our result which is \emph{simple}
enough to present a self-contained approach and \emph{general}
enough to serve as a guide for the reader in the forthcoming
papers. More precisely, we shall assume that $p=1$ and
$L_q(\mathcal{M})$ is either $\ell_q$ or the Schatten class $S_q$.
This reduces our result to the operator space embedding of
\emph{discrete} $L_q$ spaces into the predual of a von Neumann
algebra. The complete proof (including the properties of
$\mathcal{A}$) will be given in \cite{JP4}, except for the lack of
semifiniteness announced in iii) which will be considered
independently in \cite{JP5} since requires different techniques
from \cite{JP3}.

\vskip5pt

In order to put our result in the right context, let us stress the
interaction between harmonic analysis, probability and Banach
space theory carried out mostly in the 70's. Based on previous
results by Beck, Grothendieck, L{\'e}vy, Orlicz... probabilistic
methods in Banach spaces became the heart of the work developed by
Kwapie\'n, Maurey, Pisier, Rosenthal and many others. A
fundamental motivation for this new field relied on the embedding
theory of classical $L_p$ spaces. This theory was born in 1966
with the seminal paper \cite{BDK} of Bretagnolle, Dacunha-Castelle
and Krivine. They constructed an isometric embedding of $L_q$ into
$L_p$ for $1 \le p < q \le 2$, a Banach space version of our main
result; see also \cite{J0} for the analogous result with
noncommutative $L_p$ spaces. The simplest form of such embedding
was known to L\'evy and is given by
\begin{equation} \label{pstab}
\Big( \sum_{k=1}^\infty |\alpha_k|^q \Big)^{\frac1q} = \Big\|
\sum_{k=1}^\infty \alpha_k \, \theta_k \Big\|_{L_1(\Omega)},
\end{equation}
for scalars $(\alpha_k)_{k \ge 1}$ and where $(\theta_k)_{k \ge
1}$ is a suitable sequence of independent $q$-stable random
variables in $L_1(\Omega)$ for some probability space $(\Omega,
\mu)$. In other words, we have the relation $${\mathbb E} \exp
\big( i \summ_k \alpha_k \theta_k \big) = \exp \big( -c_q \summ_k
|\alpha_k|^q \big).$$

\vskip3pt

More recently, it has been discovered a parallel connection
between operator space theory and quantum probability. The
operator space version of Grothendieck theorem by Pisier and
Shlyakhtenko \cite{PS} and the embedding of $\mathrm{OH}$
\cite{J2} require tools from free probability. In this context we
should replace the $\theta_k$'s by suitable operators so that
\eqref{pstab} holds with matrix-valued coefficients $\alpha_1,
\alpha_2, \ldots$. To that aim, we develop new tools in quantum
probability and construct an operator space version of $q$-stable
random variables. To formulate this quantized form of
\eqref{pstab} we need some basic results of Pisier's theory
\cite{P2}. The most natural operator space structure on
$\ell_\infty$ comes from the diagonal embedding $\ell_\infty
\hookrightarrow \mathcal{B}(\ell_2)$. The natural structure on
$\ell_1$ is given by operator space duality, while the spaces
$\ell_p$ are defined by means of the complex interpolation method
\cite{P1} for operator spaces. Let us denote by $(\delta_k)_{k \ge
1}$ the unit vector basis of $\ell_q$. If $\widehat{\ten}$ denotes
the operator space projective tensor product and $S_p$ stands for
the Schatten $p$-class over $\ell_2$, it was shown in \cite{P2}
that
\[ \Big\| \sum_{k=1}^\infty a_k \ten \delta_k \Big\|_{S_1
\widehat{\ten} \, \ell_q} = \inf_{a_k = \alpha b_k \beta}
\|\alpha\|_{S_{2q'}} \Big( \sum_{k=1}^\infty \|b_k\|_{S_q}^q
\Big)^{\frac1q} \|\beta\|_{S_{2q'}}. \] The answer to Pisier's
problem for $\ell_q$ reads as follows.

\begin{theoA} \label{lp}
If $1 < q \le 2$, there exists a sufficiently large von Neumann
algebra $\mathcal{A}$ and a sequence $(x_k)_{k \ge 1}$ in
$L_1(\mathcal{A})$ such that the equivalence below holds for any
family $(a_k)_{k \ge 1}$ of trace class operators $$\inf_{a_k =
\alpha b_k \beta} \|\alpha\|_{2q'} \Big( \sum_{k=1}^\infty
\|b_k\|_q^q \Big)^{\frac1q} \|\beta\|_{2q'} \sim_c \Big\|
\sum_{k=1}^\infty a_k \ten x_k \Big\|_{L_1( \mathcal{A} \bar\ten
\mathcal{B}(\ell_2))}.$$
\end{theoA}

This gives a completely isomorphic embedding of $\ell_q$ into
$L_1(\mathcal{A})$. Moreover, the sequence $x_1, x_2, \ldots$
provides an operator space version of a $q$-stable process and
motivates a cb-embedding theory of $L_p$ spaces. A particular case
of Theorem A is the recent construction \cite{J2} of a
cb-embedding of Pisier's operator Hilbert space $\mathrm{OH}$ into
a von Neumann algebra predual. In other words, a complete
embedding of $\ell_2$ with its natural operator space structure
into a noncommutative $L_1$ space, see also Pisier's paper
\cite{P4} for a shorter proof and Xu's alternative construction
\cite{X2}. Other related results appear in \cite{J3,P5,PS,X3},
while semi-complete embeddings between vector-valued $L_p$ spaces
can be found in \cite{JP,Pa}. All these papers will play a role in
our analysis, either here or in the forthcoming papers.

\vskip5pt

The construction which leads to this operator space version of
$q$-stable random variables is simpler than the one needed for the
general case and will serve as a model for the latter in
\cite{JP4}. Let us sketch it in some detail. A key ingredient in
our proof is the notion of the Haagerup tensor product $\ten_h$.
We first note that $\ell_q$ is the diagonal subspace of the
Schatten class $S_q$. According to \cite{P2}, $S_q$ can be written
as the Haagerup tensor product of its first column and first row
subspaces $S_q = C_q \otimes_h R_q$. Moreover, using a simple
generalization of \lq\lq Pisier's exercise\rq\rq${}$ (see Exercise
7.9 in Pisier's book \cite{P3}) we have
\begin{equation} \label{Pexercise}
\begin{array}{rcl} C_q & \hookrightarrow_{cb} & \big( R \oplus
\mathrm{OH} \big) \big/ graph(\Lambda_1)^\perp, \\ [5pt] R_q &
\hookrightarrow_{cb} & \big( C \oplus \mathrm{OH} \big) \big/
graph(\Lambda_2)^\perp,
\end{array}
\end{equation}
with $\Lambda_1: C \to \mathrm{OH}$ and $\Lambda_2: R \to
\mathrm{OH}$ suitable injective, closed, densely-defined operators
with dense range and where $\hookrightarrow_{cb}$ denotes a
cb-embedding, see \cite{J2,P4,X3} for related results. By
\cite{J2} and duality, it suffices to see that $$graph(\Lambda_1)
\otimes_h graph(\Lambda_2)$$ is cb-isomorphic to a cb-complemented
subspace of $L_\infty(\mathcal{A};\mathrm{OH})$. Before proceeding
it is important to have a little discussion on this space. Namely,
our methods in this paper will lead us to obtain such cb-embedding
for a free product von Neumann algebra $\mathcal{A}$. In
particular, $\mathcal{A}$ will not be hyperfinite and Pisier's
theory \cite{P2} of vector valued noncommutative $L_p$ spaces does
not consider this space. Our definition of it is given by complex
interpolation $$L_\infty(\mathcal{A}; \mathrm{OH}_n) = \big[
C_n(\mathcal{A}), R_n(\mathcal{A}) \big]_{\frac12},$$ where
$R_n(\mathcal{A})$ and $C_n(\mathcal{A})$ are the row and column
subspaces of $$\mathrm{M}_n(\mathcal{A}) = \mathrm{M}_n
\otimes_{\mathrm{min}} \mathcal{A} = C_n \otimes_h \mathcal{A}
\otimes_h R_n.$$ In fact, given $1 \le p \le \infty$ and using
instead the row and column subspaces of the vector-valued Schatten
class $S_p^n(L_p(\mathcal{A}))$, we find a definition of
$L_p(\mathcal{A}; \mathrm{OH}_n)$. On the other hand, we also know
from \cite{J1} a definition of the spaces $L_p(\mathcal{A};
\ell_1^n)$ and $L_p(\mathcal{A}; \ell_\infty^n)$ for any von
Neumann algebra $\mathcal{A}$. In particular, we might wonder
whether or not our definition of $L_p(\mathcal{A}; \mathrm{OH}_n)$
satisfies the following complete isometry
$$L_p(\mathcal{A}; \mathrm{OH}_n) = \big[
L_p(\mathcal{A}; \ell_\infty^n), L_p(\mathcal{A}; \ell_1^n)
\big]_{\frac12}.$$ Fortunately this is the case and the proof can
be found in Chapter 1 of \cite{JP2}.

\vskip5pt

Let us go on with the argument. By the injectivity of the Haagerup
tensor product, $graph(\Lambda_1) \otimes_h graph(\Lambda_2)$ is
an intersection of four spaces. Let us explain this in detail. By
a standard discretization argument which will be given in Lemma
\ref{Lemma-Diagonalization} below, we may assume that $\Lambda_j =
\mathsf{d}_{\lambda^j} = \sum_k \lambda_k^j e_{kk}$ is a diagonal
operator on $\ell_2$ for $j=1,2$. Moreover, letting $[w]$ denote
the integer part of $w$ and taking
$$\lambda_k = \lambda_{[k+1/2]}^j \quad  \mbox{for} \quad k \equiv j \
(\mathrm{mod} \, 2),$$ it is no restriction to assume that the
eigenvalues of $\Lambda_1$ and $\Lambda_2$ are the same. In
particular, our considerations allow us to rewrite
$graph(\Lambda_1) \otimes_h graph(\Lambda_2)$ in the form below
$$\mathcal{J}_{\infty,2} = graph(\mathsf{d}_\lambda) \otimes_h
graph(\mathsf{d}_\lambda) = \Big( C \cap \ell_2^{oh}(\lambda)
\Big) \otimes_h \Big( R \cap \ell_2^{oh}(\lambda) \Big),$$ where
$\ell_2^{oh}(\lambda)$ is a weighted form of $\mathrm{OH}$
according to the action of $\mathsf{d}_\lambda$, so that
\begin{eqnarray*}
C \cap \ell_2^{oh} (\lambda) & = & \Big\{ ( \hskip0.5pt e_{i1}
\hskip0.5pt, \lambda_i \hskip0.5pt e_{i1} \hskip1pt) \,
\hskip0.5pt \big| \ \hskip0.5pt i \hskip0.5pt \ge 1 \Big\} \subset
C \oplus \mathrm{OH}, \\ R \cap \ell_2^{oh} (\lambda) & = & \Big\{
(e_{1j}, \lambda_j e_{1j}) \, \big| \ j \ge 1 \Big\} \subset R
\hskip0.5pt \oplus \mathrm{OH}.
\end{eqnarray*}
The notation for $\mathcal{J}_{\infty,2}$ follows \cite{JP2}.
Namely, the symbol $\infty$ in $\mathcal{J}_{\infty,2}$ is used
because we shall consider $L_p$ versions of these spaces in
\cite{JP4}. The number $2$ denotes that this space arises as a \lq
middle point\rq${}$ in the sense of interpolation theory between
two related $\mathcal{J}$-spaces, see \cite{JP2} and \cite{JP4}
for more details. Now, regarding $\mathsf{d}_\lambda^4 = \sum_k
\lambda_k^4 e_{kk}$ as the density $d_\psi$ of some normal
strictly semifinite faithful weight $\psi$ on
$\mathcal{B}(\ell_2)$, the space $\mathcal{J}_{\infty,2}$ splits
up into a $4$-term intersection space. In other words, we find
$$\mathcal{J}_{\infty,2}(\psi) = \big( C \ten_h R \big) \cap
\big( C \ten_h \mathrm{OH} \big) d_\psi^{\frac14} \cap
d_\psi^{\frac14} \big( \mathrm{OH} \ten_h R \big) \cap
d_\psi^{\frac14} \big( \mathrm{OH} \ten_h \mathrm{OH} \big)
d_\psi^{\frac14}.$$ The norm of $x$ in
$\mathcal{J}_{\infty,2}(\psi)$ is given by
$$\max \Big\{ \|x\|_{B(\ell_2)}, \big\|x d_\psi^{\frac14}
\big\|_{C \ten_h \mathrm{OH}}, \big\| d_\psi^{\frac14} x
\big\|_{\mathrm{OH} \ten_h R}, \big\| d_\psi^{\frac14} x
d_\psi^{\frac14} \big\|_{\mathrm{OH} \ten_h \mathrm{OH}} \Big\}.$$
The two middle terms are not as unusual as it might seem
\begin{equation} \label{osstruct}
\begin{array}{rcl}
\big\| d_\psi^{\frac14} (x_{ij}) \big\|_{\mathrm{M}_m(\mathrm{OH}
\ten_h R)} & = & \displaystyle \sup_{\|\alpha\|_{S_4^m}\le 1}
\Big\| d_\psi^{\frac14} \Big( \sum_{k=1}^m \alpha_{ik} x_{kj}
\Big) \Big\|_{L_4(\mathrm{M}_m \ten \mathcal{B}(\ell_2))}, \\
\big\| (x_{ij}) d_\psi^{\frac14} \big\|_{\mathrm{M}_m(C \ten_h
\mathrm{OH})} & = & \displaystyle \sup_{\|\beta\|_{S_4^m}\le 1}
\Big\| \Big( \sum_{k=1}^m x_{ik} \beta_{kj} \Big) \hskip1pt
d_\psi^{\frac14} \Big\|_{L_4(\mathrm{M}_m \ten
\mathcal{B}(\ell_2))}.
\end{array}
\end{equation}
Let us now assume that we just try to embed the finite-dimensional
space $S_q^m$. By approximation, it suffices to consider only
finitely many eigenvalues $\lambda_1, \lambda_2, \ldots,
\lambda_n$ and according to the results from \cite{J3}, we can
take $n \sim m \log m$. In this case we rename $\psi$ by $\psi_n$
and we may easily assume that $$\mbox{tr} (d_{\psi_n}) = \summ_k
\lambda_k^4 = \mathrm{k}_n$$ is a positive integer. Therefore, we
consider the following state on $\mathcal{B}(\ell_2^n)$
$$\varphi_n(x) = \frac{1}{\mathrm{k}_n} \sum_{k=1}^n \lambda_k^4
\, x_{kk}.$$ In this particular case, the space
$\mathcal{J}_{\infty,2}(\psi_n)$ can be obtained using free
probability.

\begin{theoB} \label{MAMSth}
Let $\mathcal{A} = \mathsf{A}_1 * \mathsf{A}_2 * \cdots *
\mathsf{A}_{\mathrm{k}_n}$ be the reduced free product of
$\mathrm{k}_n$ copies of $\mathcal{B}(\ell_2^n) \oplus
\mathcal{B}(\ell_2^n)$ equipped with the state $\frac12 (\varphi_n
\oplus \varphi_n)$. If $\pi_k: \mathsf{A}_j \to \mathcal{A}$
denotes the canonical embedding into the $j$-th component of
$\mathcal{A}$, the mapping $$u_n: x \in
\mathcal{J}_{\infty,2}(\psi_n) \mapsto \sum_{j=1}^{\mathrm{k}_n}
\pi_j(x,-x) \ten \delta_j \in L_\infty(\mathcal{A};
\mathrm{OH}_{\mathrm{k}_n})$$ is a cb-embedding with
cb-complemented image and constants independent of $n$.
\end{theoB}

Theorem B and its generalization for arbitrary von Neumann
algebras is a very recent result from \cite{JP2}. However, since
the result proved there covers a wider range of indices, the proof
is rather long and quite technical. In order to be self-contained
we provide a second proof of this particular case only using
elementary tools from free probability. We think this approach is
of independent interest. Now, by duality we obtain a cb-embedding
of $S_q^m$ into $L_1(\mathcal{A}; \mathrm{OH}_{\mathrm{k}_n})$.
Then, we conclude using the cb-embedding of $\mathrm{OH}$ from
\cite{J2} and an ultraproduct procedure. In fact, what we shall
prove is a far reaching generalization of Theorem \ref{lp}.
Namely, the same result holds replacing $C_q$ and $R_q$ by
subspaces of quotients of $R \oplus \mathrm{OH}$ and $C \oplus
\mathrm{OH}$.

\begin{theoC} \label{QSQSQS}
Let $\mathrm{X}_1$ be a subspace of a quotient of $R \oplus
\mathrm{OH}$ and let $\mathrm{X}_2$ be a subspace of a quotient of
$C \oplus \mathrm{OH}$. Then, there exist some $\mathrm{QWEP}$ von
Neumann algebra $\mathcal{A}$ and a cb-embedding
$$\mathrm{X}_1 \otimes_h \mathrm{X}_2 \hookrightarrow_{cb}
L_1(\mathcal{A}).$$
\end{theoC}

Now we may explain how the paper is organized. In Section
\ref{NewSect1} we just prove the complete embeddings
\eqref{Pexercise}. This is a simple consequence of Pisier's
exercise 7.9 in \cite{P3} and Xu's generalization \cite{X3}.
However, we have decided to include the proof since it will serve
as a model to follow in \cite{JP4}. In Section \ref{S2Free} we
will give an elementary proof of Theorem B using just free
probability tools. The link between Theorems B and C is given in
the first half of Section \ref{S3CBE}. The second half is devoted
to the proof of Theorem C and thereby of Theorem A.

\subsection*{Background and notation}

We shall assume that the reader is familiar with those branches of
operator algebra related to the theories of operator spaces and
noncommutative $L_p$ spaces. The recent monographs \cite{ER} and
\cite{P3} on operator spaces contain more than enough information
for our purposes. We shall work over general von Neumann algebras
so that we use Haagerup's definition \cite{H} of $L_p$, see also
Terp's excellent exposition of the subject \cite{T1}. The basics
on von Neumann algebras and Tomita's modular theory to work with
these notions appear in Kadison/Ringrose books \cite{KR}. We shall
also assume certain familiarity with Pisier's vector-valued
non-commutative $L_p$ spaces \cite{P2} and Voiculescu's free
probability theory \cite{VDN}.

\vskip5pt

We shall follow the standard notation in the subject. Anyway, let
us say a few words on our terminology. The symbols $(\delta_k)$
and $(e_{ij})$ will denote the unit vector basis of $\ell_2$ and
$\mathcal{B}(\ell_2)$ respectively. The letters $\mathcal{A},
\mathcal{M}$ and $\mathcal{N}$ are reserved to denote von Neumann
algebras. Almost all the time, the inclusions $\mathcal{N} \subset
\mathcal{M} \subset \mathcal{A}$ will hold. We shall use $\varphi$
and $\phi$ to denote normal faithful (\emph{n.f.} in short)
states, while the letter $\psi$ will be reserved for normal
strictly semifinite faithful (\emph{n.s.s.f.} in short) weights.
Inner products and duality brackets will be anti-linear on the
first component and linear on the second. Given $\gamma > 0$, we
shall write $\gamma \mathrm{X}$ to denote the space $\mathrm{X}$
equipped with the norm $\|x\|_{\gamma \mathrm{X}} = \gamma
\|x\|_{\mathrm{X}}$. In particular, if $\mathcal{M}$ is a finite
von Neumann algebra equipped with a finite weight $\psi$, we shall
usually write $\psi = \mathrm{k} \varphi$ with $\mathrm{k} =
\psi(1)$ so that $\varphi$ becomes a state on $\mathcal{M}$. In
this situation, the associated $L_p$ space will be denoted as
$\mathrm{k}^{1/p} L_p(\mathcal{M})$, so that the $L_p$ norm is
calculated using the state $\varphi$. Note that this notation does
not follow the tradition in interpolation theory, which denotes by
$\gamma \mathrm{X}$ the Banach space with unit ball $\gamma
\mathsf{B}_{\mathrm{X}}$ and so the norm is divided and not
multiplied by $\gamma$. Given an operator space $\mathrm{X}$, the
$\mathrm{X}$-valued Schatten $p$-class over the algebra
$\mathcal{B}(\ell_2)$ will be denoted by $S_p(\mathrm{X})$.

\subsection*{Acknowledgement}

The former version of this paper was much longer and also
significantly more technical, since it almost contained the proof
of our main result in full generality. This included all the
results in this paper and also the results announced for
\cite{JP4}. We owe the referee to suggest this more accessible
presentation which we hope will help the nonexpert reader.

\numberwithin{equation}{section}

\section{On \lq\lq Pisier's exercise\rq\rq}
\label{NewSect1}

We begin by proving a generalization of Exercise 7.9 in \cite{P3}.
This result became popular after Pisier applied it in \cite{P4} to
obtain a simpler way to cb-embed $\mathrm{OH}$ into the predual of
a von Neumann algebra. In fact, our argument is quite close to the
one given in \cite{X3} for a similar result and might be known to
experts. Nevertheless we include it here for completeness and also
since we can think of it as a model for our construction in the
non-discrete case. In fact, in this paper we shall only use Lemma
\ref{Lemma-Motivation} below for $p=1$, while in \cite{JP4} we
will apply it in full generality.

\vskip5pt

Let us set some notation. Given a Hilbert space $\mathcal{H}$, we
shall write $\mathcal{H}_r = \mathcal{B}(\mathcal{H}, \C)$ and
$\mathcal{H}_c = \mathcal{B}(\C, \mathcal{H})$ for the row and
column quantizations on $\mathcal{H}$. Moreover, given $1 \le p
\le \infty$ we shall use the following terminology
$$\mathcal{H}_{r_p} = \big[ \mathcal{H}_r, \mathcal{H}_c
\big]_{\frac1p} \quad \mbox{and} \quad \mathcal{H}_{c_p} = \big[
\mathcal{H}_c, \mathcal{H}_r\big]_{\frac1p}.$$ There are two
particular cases for which we use another terminology. When
$\mathcal{H} = \ell_2$, we shall use $(R,C,R_p,C_p)$ instead.
Moreover, when the Hilbert space is $L_2(\mathcal{M})$ for some
von Neumann algebra $\mathcal{M}$, we shall write
$L_2^{r_p}(\mathcal{M})$ and $L_2^{c_p}(\mathcal{M})$. In the same
fashion, $\mathcal{H}_{oh}$ and $L_2^{oh}(\mathcal{M})$ stand for
the operator Hilbert space structures. Given two operator spaces
$\mathrm{X}_1$ and $\mathrm{X}_2$, the expression $\mathrm{X}_1
\simeq_{cb} \mathrm{X}_2$ means that there exists a complete
isomorphism between them. We shall write $\mathrm{X}_1 \in
\mathcal{QS}(\mathrm{X}_2)$ to denote that $\mathrm{X}_1$ is
completely isomorphic to a quotient of a subspace of
$\mathrm{X}_2$. Let $\mathcal{S}$ denote the strip
$$\mathcal{S} = \Big\{ z \in \C \, \big| \ 0 < \mbox{Re}(z) < 1
\Big\}$$ and let $\partial \mathcal{S} =
\partial_0 \cup \partial_1$ be the partition of its boundary into
$$\partial_0 = \Big\{ z \in \C \, \big| \ \mbox{Re}(z)=0 \Big\}
\quad \mbox{and} \quad \partial_1 = \Big\{ z \in \C \, \big| \
\mbox{Re}(z)=1 \Big\}.$$ Given $0 < \theta < 1$, let $\mu_\theta$
be the harmonic measure of the point $z = \theta$. This is a
probability measure on $\partial \mathcal{S}$ (with density given
by the Poisson kernel in the strip) that can be written as
$\mu_\theta = (1 - \theta) \mu_0 + \theta \mu_1$, with $\mu_j$
being probability measures supported by $\partial_j$ and such that
\begin{equation} \label{Eq-AnalyticCondition}
f(\theta) = \int_{\partial \mathcal{S}} f d\mu_\theta
\end{equation}
for any bounded analytic function $f: \mathcal{S} \to \C$ which is
extended non-tangentially to $\partial \mathcal{S}$. Now, before
proving the announced result, we need to set a formula from
\cite{X}. Given $2 \le p \le \infty$ and $0 < \theta < 1$, the
norm of $x = \sum_k x_k \ten \delta_k$ in $$\big[ S_p(C_p),
S_p(R_p) \big]_{\theta}$$ is given by
\begin{equation} \label{fbasexu}
\sup \Big\{ \Big( \summ_k \big\| \alpha x_k \beta \big\|_{S_2}^2
\Big)^{\frac12} \big| \ \|\alpha\|_{S_u}, \|\beta\|_{S_v} \le 1
\Big\},
\end{equation}
where the indices $(u,v)$ are determined by $$(1/u,1/v) =
(\theta/q, (1-\theta)/q) \quad \mbox{with} \quad 1/2 = 1/p +
1/q.$$ Of course, this formula trivially generalizes for the norm
in the space $$\big[ S_p(\mathcal{H}_{c_p}),
S_p(\mathcal{H}_{r_p}) \big]_{\theta}.$$

\begin{lemma} \label{Lemma-Motivation}
If $1 \le p < q \le 2$, we have
$$R_q \in \mathcal{QS} \big( R_p \oplus_2
\mathrm{OH} \big) \quad \mbox{and} \quad C_q \in \mathcal{QS}
\big( C_p \oplus_2 \mathrm{OH} \big).$$
\end{lemma}

\dem We only prove the first assertion since the arguments for
both are the same. In what follows we fix $0 < \theta < 1$
determined by the relation $R_q = [R_p, \mathrm{OH}]_\theta$. In
other words, we have $1/q = (1-\theta)/p + \theta/2$. According to
the complex interpolation method and its operator space extension
\cite{P1}, given a compatible couple $(\mathrm{X}_0,
\mathrm{X}_1)$ of operator spaces we define
$\mathcal{F}(\mathrm{X}_0, \mathrm{X}_1)$ as the space of bounded
analytic functions $f: \mathcal{S} \to \mathrm{X}_0 +
\mathrm{X}_1$ and we equip it with the following norm
$$\|f\|_{\mathcal{F}(\mathrm{X}_0, \mathrm{X}_1)} = \Big(
(1-\theta) \,
\|{f_{\mid_{\partial_0}}}\|_{L_2(\partial_0;\mathrm{X}_0)}^2 +
\theta \,
\|{f_{\mid_{\partial_1}}}\|_{L_2(\partial_1;\mathrm{X}_1)}^2
\Big)^{\frac12}.$$ Then, the complex interpolation space
$\mathrm{X}_\theta = [\mathrm{X}_0, \mathrm{X}_1]_{\theta}$ can be
defined as the space of all $x \in \mathrm{X}_0 + \mathrm{X}_1$
such that there exists a function $f \in \mathcal{F}(\mathrm{X}_0,
\mathrm{X}_1)$ with $f(\theta) = x$. We equip $\mathrm{X}_\theta$
with the norm $$\|x\|_{\mathrm{X}_\theta} = \inf \Big\{
\|f\|_{\mathcal{F}(\mathrm{X}_0, \mathrm{X}_1)} \, \big| \ f \in
\mathcal{F}(\mathrm{X}_0, \mathrm{X}_1) \ \mbox{and} \ f(\theta) =
x \Big\}.$$ In our case we set $\mathrm{X}_0 = R_p$ and
$\mathrm{X}_1 = \mathrm{OH}$. If we define $$\mathcal{H} =
(1-\theta)^{\frac12} L_2(\partial_0; \ell_2) \quad \mbox{and}
\quad \mathcal{K} = \theta^\frac12 L_2(\partial_1; \ell_2),$$ it
turns out that $\mathcal{F}(R_p,\mathrm{OH})$ can be regarded (via
Poisson integration) as a subspace of $\mathcal{H} \oplus_2
\mathcal{K}$. Moreover, we equip $\mathcal{F}(R_p,\mathrm{OH})$
with the operator space structure inherited from
$\mathcal{H}_{r_p} \oplus_2 \mathcal{K}_{oh}$. Then, we define the
mapping $$\mathcal{Q}: f \in \mathcal{F}(R_p,\mathrm{OH}) \mapsto
f(\theta) \in R_q.$$ The assertion follows from the fact that
$\mathcal{Q}$ is a complete metric surjection. Indeed, in that
case we have $R_q \simeq_{cb} \mathcal{F}(R_p,\mathrm{OH}) / \ker
\mathcal{Q}$, which is a quotient of a subspace of $R_p \oplus_2
\mathrm{OH}$. In order to see that $\mathcal{Q}$ is a complete
surjection, it suffices to see that the map $id_{S_{p'}} \otimes
\mathcal{Q}: S_{p'}(\mathcal{F}(R_p,\mathrm{OH})) \to S_{p'}(R_q)$
is a metric surjection. We begin by showing that $id_{S_{p'}}
\otimes \mathcal{Q}$ is contractive. Let $f \in
S_{p'}(\mathcal{F}(R_p,\mathrm{OH}))$ be of norm $< 1$ and let us
write $$f(\theta) = \summ_k f_k(\theta) \otimes \delta_k \in
S_{p'}(R_q).$$ To compute the norm of $f(\theta)$ we note that
$$S_{p'}(R_q) = \big[ S_{p'}(C_{p'}), S_{p'}(R_{p'}) \big]_\eta
\quad \mbox{with} \quad 1/q = (1-\eta)/p + \eta/p'.$$ Then it
follows from \eqref{fbasexu} that
\begin{equation} \label{Eq-int-f(theta)}
\|f(\theta)\|_{S_{p'}(R_q)} = \sup \Big\{ \Big( \summ_k \big\|
\alpha f_k(\theta) \beta \big\|_{S_2}^2 \Big)^{\frac12} \, \big| \
\|\alpha\|_{S_{2r/\eta}}, \|\beta\|_{S_{2r/(1-\eta)}} \le 1 \Big\}
\end{equation}
where $1/2r = 1/2 - 1/p' = 1/p - 1/2$. Moreover, by polar
decomposition it is clear that we can restrict the supremum above
to all $\alpha$ and $\beta$ in the positive parts of their
respective unit balls. Taking this restriction in consideration,
we define
$$g(z) = \summ_k g_k(z) \otimes \delta_k \quad \mbox{with} \quad
g_k(z) = \alpha^{\frac{z}{\theta}} f_k(z)
\beta^{\frac{2-z}{2-\theta}}.$$ The $g_k$'s are analytic in
$\mathcal{S}$ and take values in $S_2$. Thus, we have the identity
\begin{eqnarray} \label{Eq-id-harm}
\lefteqn{\Big( \summ_k \big\| \alpha f_k(\theta) \beta
\big\|_{S_2}^2 \Big)^{\frac12}} \\ \nonumber & = & \Big( \summ_k
\|g_k(\theta)\|_{S_2}^2 \Big)^{\frac12} \\ \nonumber & = & \Big(
(1-\theta) \int_{\partial_0} \summ_k \|g_k(z)\|_{S_2}^2 d\mu_0 +
\theta \int_{\partial_1} \summ_k \|g_k(z)\|_{S_2}^2 d\mu_1
\Big)^{\frac12}.
\end{eqnarray}
The contractivity of $id_{S_{p'}} \otimes \mathcal{Q}$ will follow
from
\begin{eqnarray}
\label{Est000} \int_{\partial_0} \summ_k \|g_k(z)\|_{S_2}^2 d\mu_0
& \le & \|f_{\mid_{\partial_0}}\|_{S_{p'}
(L_2^{r_p}(\partial_0;\ell_2))}^2, \\ \label{Est111}
\int_{\partial_1} \summ_k \|g_k(z)\|_{S_2}^2 d\mu_1 & \le &
\|f_{\mid_{\partial_1}}\|_{S_{p'}(L_2^{oh}(\partial_1;\ell_2))}^2.
\end{eqnarray}
Indeed, if we combine \eqref{Eq-int-f(theta)} and
\eqref{Eq-id-harm} with the operator space structure defined on
$\mathcal{F}(R_p, \mathrm{OH})$, it turns out that inequalities
\eqref{Est000} and \eqref{Est111} are exactly what we need. To
prove \eqref{Est000} we observe that $2 \eta = \theta$ and
$r'=p'/2$, so that
\begin{eqnarray*}
\int_{\partial_0} \summ_k \|g_k(z)\|_{S_2}^2 d\mu_0 & = &
\int_{\partial_0} \summ_k \|f_k(z)
\beta^{\frac{2}{2-\theta}}\|_{S_2}^2 d\mu_0 \\ & = &
\int_{\partial_0} \summ_k \mbox{tr} \big( f_k(z)^* f_k(z)
\beta^{\frac{1}{1-\eta}} (\beta^{\frac{1}{1-\eta}})^* \big) d\mu_0
\\ & \le & \|\beta^{\frac{2}{1-\eta}}\|_{S_r} \Big\|
\int_{\partial_0} \summ_k f_k(z)^* f_k(z) d\mu_0 \Big\|_{S_{r'}}.
\end{eqnarray*}
This gives $$\int_{\partial_0} \summ_k \|g_k(z)\|_{S_2}^2 d\mu_0
\le \|\beta\|_{S_{2r/(1-\eta)}}^{\frac{2}{1-\eta}} \Big\| \Big(
\int_{\partial_0} \summ_k f_k(z)^* f_k(z) d\mu_0 \Big)^{\frac12}
\Big\|_{S_{p'}}^2.$$ The first term on the right is $\le 1$.
Hence, we have
\begin{eqnarray*}
\int_{\partial_0} \summ_k \|g_k(z)\|_{S_2}^2 d\mu_0 & \le &
\|f_{\mid_{\partial_0}}\|_{S_{p'}(L_2^{c_{p'}}(\partial_0;
\ell_2))}^2 = \|f_{\mid_{\partial_0}}\|_{S_{p'}
(L_2^{r_p}(\partial_0;\ell_2))}^2.
\end{eqnarray*}
This proves \eqref{Est000} while for \eqref{Est111} we proceed in
a similar way
\begin{eqnarray*}
\lefteqn{\int_{\partial_1} \summ_k \|g_k(z)\|_{S_2}^2 d\mu_1} \\ &
= & \int_{\partial_1} \summ_k \| \alpha^{\frac{1}{\theta}} f_k(z)
\beta^{\frac{1}{2-\theta}}\|_{S_2}^2 d\mu_1 \\ & = &
\int_{\partial_1} \summ_k \| \alpha^{\frac{1}{2 \eta}} f_k(z)
\beta^{\frac{1}{2- 2 \eta}}\|_{S_2}^2 d\mu_1 \\ & \le & \sup
\left\{ \int_{\partial_1} \summ_k \| a f_k(z) b \|_{S_2}^2 d\mu_1
\, \big| \ \|a\|_{S_{4r}}, \|b\|_{S_{4r}} \le 1 \right\}
\\ & = & \sup \left\{ \Big\| a \Big(
\summ_k {f_k}_{\mid_{\partial_1}} \otimes \delta_k \Big) b
\Big\|_{S_2(L_2^{oh}(\partial_1; \ell_2))}^2 \, \big| \
\|a\|_{S_{4r}}, \|b\|_{S_{4r}} \le 1 \right\}.
\end{eqnarray*}
According to \eqref{fbasexu}, this proves \eqref{Est111} and we
have a contraction. Reciprocally, given $x \in S_{p'}(R_q)$ of
norm $< 1$, we are now interested on finding $f \in
S_{p'}(\mathcal{F}(R_p, \mathrm{OH}))$ such that $f(\theta) = x$
and $\|f\|_{S_{p'}(\mathcal{F}(R_p, \mathrm{OH}))} \le 1$. Since
$$[S_{p'}(R_p), S_{p'}(\mathrm{OH})]_\theta = S_{p'}(R_q)$$ there
must exists $f \in \mathcal{F}(S_{p'}(R_p), S_{p'}(\mathrm{OH}))$
such that $f(\theta) = x$ and
$$\|f\|_{\mathcal{F}(S_{p'}(R_p),
S_{p'}(\mathrm{OH}))} = \Big( (1-\theta)
\|f_{\mid_{\partial_0}}\|_{L_2(\partial_0;S_{p'}(R_p))}^2 + \theta
\|f_{\mid_{\partial_1}}\|_{L_2(\partial_1;S_{p'}(\mathrm{OH}))}^2
\Big)^{\frac12} \le 1.$$ Therefore, it remains to see that
\begin{eqnarray}
\label{Eq-222}
\|f_{\mid_{\partial_0}}\|_{S_{p'}(L_2^{r_p}(\partial_0;\ell_2))} &
\le & \|f_{\mid_{\partial_0}}\|_{L_2(\partial_0;S_{p'}(R_p))}, \\
\label{Eq-333}
\|f_{\mid_{\partial_1}}\|_{S_{p'}(L_2^{oh}(\partial_1;\ell_2))} &
\le &
\|f_{\mid_{\partial_1}}\|_{L_2(\partial_1;S_{p'}(\mathrm{OH}))}.
\end{eqnarray}
However, the identities below are clear by now
\begin{eqnarray*}
\|f_{\mid_{\partial_0}}\|_{S_{p'}(L_2^{r_p}(\partial_0;\ell_2))} &
= & \Big\| \Big( \int_{\partial_0} \summ_k f_k(z)^* f_k(z) d\mu_0
\Big)^{\frac12}\Big\|_{S_{p'}}, \\
\|f_{\mid_{\partial_0}}\|_{L_2 (\partial_0; S_{p'}(R_p))}
\hskip1.5pt & = & \Big( \int_{\partial_0} \Big\| \Big( \summ_k
f_k(z)^* f_k(z) \Big)^{\frac12} \Big\|_{S_{p'}}^2 d\mu_0
\Big)^{\frac12}.
\end{eqnarray*}
In particular, \eqref{Eq-222} follows automatically. On the other
hand, \eqref{fbasexu} gives
\begin{eqnarray*}
\|f_{\mid_{\partial_1}}\|_{S_{p'}(L_2^{oh}(\partial_1;\ell_2))} &
= & \sup_{\|\alpha\|_{4r}, \|\beta\|_{4r} \le 1} \Big(
\int_{\partial_1} \summ_k \|\alpha f_k(z) \beta\|_{S_2}^2 d\mu_1
\Big)^{\frac12}, \\
\|f_{\mid_{\partial_1}}\|_{L_2 (\partial_1; S_{p'}(\mathrm{OH}))}
& = & \Big( \int_{\partial_1} \sup_{\|\alpha\|_{4r},
\|\beta\|_{4r} \le 1} \summ_k \|\alpha f_k(z) \beta\|_{S_2}^2
d\mu_1 \Big)^{\frac12}.
\end{eqnarray*}
Thus, inequality \eqref{Eq-333} also follows easily and
$\mathcal{Q}$ is a complete metric surjection. \fin

\begin{remark}
\emph{If $1 \le p_0 \le p \le p_1 \le \infty$, it is also true
that $$R_p \in \mathcal{QS} \big( R_{p_0} \oplus_2 R_{p_1} \big)
\quad \mbox{and} \quad C_p \in \mathcal{QS} \big( C_{p_0} \oplus_2
C_{p_1} \big).$$}
\end{remark}

\begin{remark} \label{Remark-Graph_Sq}
\emph{In Lemma \ref{Lemma-Motivation} we have obtained $$R_q
\simeq_{cb} \mathcal{F}(R_p, \mathrm{OH}) / \ker \mathcal{Q} \in
\mathcal{QS}(R_p \oplus_2 \mathrm{OH}),$$ $$C_q \simeq_{cb}
\mathcal{F}(C_p, \mathrm{OH}) / \ker \mathcal{Q} \in
\mathcal{QS}(C_p \oplus_2 \mathrm{OH}).$$ However, it will be
convenient in the sequel to observe that $\ker \mathcal{Q}$ can be
regarded in both cases as the annihilator of the graph of certain
linear operator, see \eqref{Pexercise}. Recall that for a linear
map between Hilbert spaces $\Lambda: \mathcal{K}_1 \to
\mathcal{K}_2$ with domain $\mathsf{dom}(\Lambda)$
$$graph (\Lambda) = \Big\{ (x_1,x_2) \in \mathcal{K}_1 \oplus_2
\mathcal{K}_2 \, \big| \ x_1 \in \mathsf{dom}(\Lambda) \
\mbox{and} \ x_2 = \Lambda(x_1) \Big\}.$$ Let us consider for
instance the case of $R_q$. We first note that
$\mathcal{F}(R_p,\mathrm{OH})$ is the graph of an injective,
closed, densely-defined operator $\Lambda$ with dense range. This
operator is given by $\Lambda(f_{\mid_{\partial_0}}) =
f_{\mid_{\partial_1}}$ for all $f \in \mathcal{F}(R_p,
\mathrm{OH})$. $\Lambda$ is well defined and injective by Poisson
integration due to the analyticity of elements in
$\mathcal{F}(R_p, \mathrm{OH})$. On the other hand, $\ker
\mathcal{Q}$ is the subspace of $\mathcal{F}(R_p, \mathrm{OH})$
composed of functions $f$ vanishing at $z = \theta$. Then, it
easily follows from \eqref{Eq-AnalyticCondition} that $\ker
\mathcal{Q}$ is the annihilator of $\mathcal{F}(R_{p'},
\mathrm{OH}) = graph (\Lambda)$ regarded as a subspace of
$$(1-\theta)^{\frac12} L_2^{r_{p'}}(\partial_0; \ell_2) \oplus_2
\theta^{\frac12} L_2^{oh} (\partial_1; \ell_2).$$}
\end{remark}

\section{Free harmonic analysis}
\label{S2Free}

We start with the proof of a generalized form of Theorem
\ref{MAMSth} for arbitrary von Neumann algebras. Although we just
use in this paper its discrete version as stated in the
Introduction, the general formulation does not present extra
difficulties and will be instrumental when dealing with
non-discrete algebras in \cite{JP4}. Our starting point is a von
Neumann algebra $\mathcal{M}$ equipped with a \emph{n.f.} state
$\varphi$ and associated density $d_{\varphi}$. Let $\mathcal{N}$
be a von Neumann subalgebra of $\mathcal{M}$. According to
Takesaki \cite{Ta2}, the existence and uniqueness of a \emph{n.f.}
conditional expectation $\mathsf{E}: \mathcal{M} \rightarrow
\mathcal{N}$ is equivalent to the invariance of $\mathcal{N}$
under the action of the modular automorphism group
$\sigma_t^{\varphi}$ associated to $(\mathcal{M},\varphi)$. In
that case, $\mathsf{E}$ is $\varphi$-invariant and following
Connes \cite{Co} we have $\mathsf{E} \circ \sigma_t^\varphi =
\sigma_t^\varphi \circ \mathsf{E}$. In what follows, we assume
these properties in all subalgebras considered. Now we set
$\mathsf{A}_k = \mathcal{M} \oplus \mathcal{M}$ and define
$\mathcal{A}$ to be the reduced amalgamated free product
$*_{\mathcal{N}} \mathsf{A}_k$ of the family $\mathsf{A}_1,
\mathsf{A}_2, \ldots, \mathsf{A}_n$ over the subalgebra
$\mathcal{N}$, where the embedding of the subalgebra $\mathcal{N}$
into $\mathcal{M} \oplus \mathcal{M}$ is given by $x \mapsto
(x,x)$. Note that our notation $*_\mathcal{N} \mathsf{A}_k$ for
reduced amalgamated free products does not make explicit the
dependence on the conditional expectations $\mathsf{E}_k:
\mathsf{A}_k \to \mathcal{N}$, given by $\mathsf{E}_k(a,b) =
\frac12 \mathsf{E}(a) + \frac12 \mathsf{E}(b)$. The following is
the operator-valued version \cite{J2,JPX} of Voiculescu inequality
\cite{V2}, for which we need to introduce the mean-zero subspaces
$$\bubl_k = \Big\{ x \in \mathsf{A}_k \, \big| \ \mathsf{E}_k(a_k)
= 0 \Big\}.$$

\begin{lemma} \label{Lemma-Voiculescu}
If $a_k \in \bubl_k$ for $1 \le k \le n$ and
$\mathsf{E}_\mathcal{N}: \mathcal{A} \to \mathcal{N}$ stands for
the conditional expectation of $\mathcal{A}$ onto $\mathcal{N}$,
the following equivalence of norms holds with constants
independent of $n$ $$\Big\| \sum_{k=1}^n a_k \Big\|_{\mathcal{A}}
\sim \max_{1 \le k \le n} \|x_k\|_{\mathsf{A}_k} + \Big\| \Big(
\sum_{k=1}^n \mathsf{E}_{\mathcal{N}}(a_k a_k^*) \Big)^{\frac12}
\Big\|_{\mathcal{N}} + \Big\| \Big( \sum_{k=1}^n
\mathsf{E}_{\mathcal{N}}(a_k^* a_k) \Big)^{\frac12}
\Big\|_{\mathcal{N}}.$$ Moreover, we also have
$$\Big\| \big( \sum_{k=1}^n a_k a_k^* \big)^{\frac12}
\Big\|_{\mathcal{A}} \sim \max_{1 \le k \le n}
\|a_k\|_{\mathsf{A}_k} + \Big\| \Big( \sum_{k=1}^n
\mathsf{E}_{\mathcal{N}}(a_k a_k^*) \Big)^{\frac12}
\Big\|_{\mathcal{N}},$$ $$\Big\| \big( \sum_{k=1}^n a_k^*a_k
\big)^{\frac12} \Big\|_{\mathcal{A}} \sim \max_{1 \le k \le n}
\|a_k\|_{\mathsf{A}_k} + \Big\| \Big( \sum_{k=1}^n
\mathsf{E}_{\mathcal{N}}(a_k^* a_k) \Big)^{\frac12}
\Big\|_{\mathcal{N}}.$$
\end{lemma}

\dem For the first inequality we refer to \cite{J2}. The others
can be proved in a similar way. Alternatively, both can be deduced
from the first one. Indeed, using the identity $$\Big\| \big(
\sum_{k=1}^n a_k a_k^* \big)^{\frac12} \Big\|_\mathcal{A} = \Big\|
\sum_{k=1}^n a_k \otimes e_{1k}
\Big\|_{\mathrm{M}_n(\mathcal{A})}$$ and recalling the isometric
isomorphism $$\mathrm{M}_n \big( *_{\mathcal{N}} \mathsf{A}_k
\big) = *_{\mathrm{M}_n(\mathcal{N})}
\mathrm{M}_n(\mathsf{A}_k),$$ we may apply Voiculescu's inequality
over the triple $$\Big( \mathrm{M}_n(\mathcal{A}),
\mathrm{M}_n(\mathsf{A}_k), \mathrm{M}_n(\mathcal{N}) \Big).$$
Taking $\widetilde{\mathsf{E}}_{\mathcal{N}} = id_{\mathrm{M}_n}
\ten \mathsf{E}_{\mathcal{N}}$, the last term disappears because
$$\Big\| \Big( \sum_{k=1}^n \widetilde{\mathsf{E}}_{\mathcal{N}}
\big( (a_k \ten e_{1k})^* (a_k \ten e_{1k}) \big) \Big)^{\frac12}
\Big\|_{\mathrm{M}_n(\mathcal{N})} \! = \! \sup_{1 \le k \le n}
\big\| \mathsf{E}_\mathcal{N}(a_k^* a_k)^{\frac12}
\big\|_{\mathcal{N}} \le \sup_{1 \le k \le n}
\|a_k\|_{\mathsf{A}_k}.$$ The third equivalence follows by taking
adjoints. The proof is complete. \fin

Let $\pi_k: \mathsf{A}_k \to \mathcal{A}$ denote the embedding of
$\mathsf{A}_k$ into the $k$-th component of $\mathcal{A}$. Given
$x \in \mathcal{M}$, we shall write $x_k$ as an abbreviation of
$\pi_k(x,-x)$. Note that $x_k$ is mean-zero. In the following we
shall use with no further comment the identities
$\mathsf{E}_{\mathcal{N}}(x_k x_k^*) = \mathsf{E}(x x^*)$ and
$\mathsf{E}_{\mathcal{N}}(x_k^* x_k) = \mathsf{E}(x^* x)$. We will
mostly work with identically distributed variables. In other
words, given $x \in \mathcal{M}$ we shall work with the sequence
$x_k = \pi_k(x,-x)$ for $1 \le k \le n$. In terms of the last
equivalences in Lemma \ref{Lemma-Voiculescu}, we may consider the
following norms
\begin{eqnarray*}
\|x\|_{\mathcal{R}_{\infty,1}^n} & = & \max \Big\{
\|x\|_{\mathcal{M}}, \sqrt{n} \, \big\| \mathsf{E}
(xx^*)^{\frac12} \big\|_{\mathcal{N}} \Big\}, \\
\|x\|_{\hskip1pt \mathcal{C}_{\infty,1}^n} \hskip2pt & = & \max
\Big\{ \|x\|_{\mathcal{M}}, \sqrt{n} \, \big\| \mathsf{E}
(x^*x)^{\frac12} \big\|_{\mathcal{N}} \Big\}.
\end{eqnarray*}
Here the letters $\mathcal{R}$ and $\mathcal{C}$ stand for row and
column according to Lemma \ref{Lemma-Voiculescu}. The symbol
$\infty$ is motivated because in \cite{JP2} and \cite{JP4} we work
with $L_p$ versions of these spaces. The number $1$ arises from
interpolation theory, because we think of these spaces as
endpoints in an interpolation scale. Finally, the norms on the
right induce to introduce the spaces $L_\infty^r(\mathcal{M},
\mathsf{E})$ and $L_\infty^c(\mathcal{M}, \mathsf{E})$ as the
closure of $\mathcal{M}$ with respect to the norms $$\big\|
\mathsf{E} (xx^*)^{\frac12} \big\|_{\mathcal{N}} \quad \mbox{and}
\quad \big\| \mathsf{E} (x^*x)^{\frac12} \big\|_{\mathcal{N}}.$$
In this way, we obtain the spaces
$$\begin{array}{rclcl} \mathcal{R}_{\infty,1}^n (\mathcal{M},
\mathsf{E}) & = &
\mathcal{M} \cap \sqrt{n} \, L_{\infty}^r(\mathcal{M}, \mathsf{E}), \\
[5pt] \mathcal{C}_{\infty,1}^n \, (\mathcal{M}, \mathsf{E}) & = &
\mathcal{M} \cap \sqrt{n} \, L_{\infty}^c(\mathcal{M},
\mathsf{E}). \end{array}$$

\begin{remark}
\emph{It is easily seen that
\begin{eqnarray*}
\big\| \mathsf{E} (xx^*)^{\frac12} \big\|_{\mathcal{N}} & = & \sup
\Big\{ \|\alpha x\|_{L_2(\mathcal{M})} \, \big| \
\|\alpha\|_{L_2(\mathcal{N})} \le 1 \Big\}, \\ \big\| \mathsf{E}
(x^*x)^{\frac12} \big\|_{\mathcal{N}} & = & \sup \Big\{ \|x
\beta\|_{L_2(\mathcal{M})} \, \big| \ \|\beta\|_{L_2(\mathcal{N})}
\le 1 \Big\}.
\end{eqnarray*}
This relation will be crucial in this paper and will be assumed in
what follows.}
\end{remark}

The state $\varphi$ induces the \emph{n.f.} state $\phi = \varphi
\circ \mathsf{E}_\mathcal{N}$ on $\mathcal{A}$. If
$\mathcal{A}_{\oplus n}$ denotes the $n$-fold direct sum
$\mathcal{A} \oplus \mathcal{A} \oplus \ldots \oplus \mathcal{A}$
(which we shall identify with $\mathcal{A} \otimes \C^n$ in the
sequel) we consider the \emph{n.f.} state $\phi_n:
\mathcal{A}_{\oplus n} \rightarrow \C$ and the conditional
expectation $\mathcal{E}_n: \mathcal{A}_{\oplus n} \rightarrow
\mathcal{A}$ given by $$\phi_n \big( \sum_{k=1}^n a_k \ten
\delta_k \big) = \frac{1}{n} \sum_{k=1}^n \phi(a_k) \qquad
\mbox{and} \qquad \mathcal{E}_n \big( \sum_{k=1}^n a_k \ten
\delta_k \big) = \frac{1}{n} \sum_{k=1}^n a_k.$$ Let us consider
the map
\begin{equation} \label{Eq-Map-u}
u: x \in \mathcal{M} \mapsto \sum_{k=1}^n x_k \ten \delta_k \in
\mathcal{A}_{\oplus n} \quad \mbox{with} \quad x_k = \pi_k(x,-x).
\end{equation}

\begin{lemma} \label{Lemma-Complemented-Isomorphism}
The mappings $$\begin{array}{rrcl} u_r: & x \in
\mathcal{R}_{\infty,1}^n (\mathcal{M}, \mathsf{E}) & \mapsto &
\displaystyle \sum_{k=1}^n x_k \ten e_{1k} \in R_n(\mathcal{A}), \\
[5pt] u_c: & x \in \mathcal{C}_{\infty,1}^n \, (\mathcal{M},
\mathsf{E}) & \mapsto & \displaystyle \sum_{k=1}^n x_k \ten e_{k1}
\in C_n(\mathcal{A}),
\end{array}$$ are isomorphisms onto complemented subspaces with
constants independent of $n$.
\end{lemma}

\dem Given $x \in \mathcal{R}_{\infty,1}^n (\mathcal{M},
\mathsf{E})$, Lemma \ref{Lemma-Voiculescu} gives
$$\big\| u_r(x) \big\|_{R_n(\mathcal{A})} = \Big\| \big(
\sum_{k=1}^n x_k x_k^* \big)^{\frac12} \Big\|_{\mathcal{A}} \sim
\max_{1 \le k \le n} \|x_k\|_{\mathsf{A}_k} + \Big\| \Big(
\sum_{k=1}^n \mathsf{E}_{\mathcal{N}} \big( x_k x_k^* \big)
\Big)^{\frac12} \Big\|_{\mathcal{N}}.$$ In other words, we have
$$\big\| u_r(x) \big\|_{R_n(\mathcal{A})} \sim
\|x\|_{\mathcal{M}} + \sqrt{n} \, \|x\|_{L_{\infty}^r(\mathcal{M},
\mathsf{E})} \sim \|x\|_{\mathcal{R}_{\infty,1}^n(\mathcal{M},
\mathsf{E})}.$$ Thus $u_r$ is an isomorphism onto its image with
constants independent of $n$. The same argument yields to the same
conclusion for $u_c$. Let $d_\varphi$ and $d_\phi$ be the
densities associated to the states $\varphi$ and $\phi$. To prove
the complementation, we consider the map
$$\omega_r: x \in L_1(\mathcal{M}) + \frac{1}{\sqrt{n}} \,
L_1^r(\mathcal{M}, \mathsf{E}) \longmapsto \frac{1}{n}
\sum_{k=1}^n x_k \otimes e_{1k} \in R_1^n(L_1(\mathcal{A})),$$
where $L_1^r(\mathcal{M}, \mathsf{E})$ is the closure of
$\mathcal{N} d_\varphi \mathcal{M}$ with respect to the norm
$\|\mathsf{E}(xx^*)^{\frac12}\|_1$. Now we use approximation and
assume that $x = \alpha d_\varphi a$ for some $(\alpha,a) \in
\mathcal{N} \times \mathcal{M}$. Then, taking $a_k = \pi_k(a,-a)$
it follows from Theorem 7.1 in \cite{JX} that
\begin{eqnarray*}
\big\| \omega_r(x) \big\|_{R_1^n(L_1(\mathcal{A}))} & = &
\frac{1}{n} \, \Big\| \alpha d_\phi \Big( \sum_{k=1}^n a_k a_k^*
\Big) d_\phi \alpha^* \Big\|_{L_{1/2}(\mathcal{A})}^{1/2}
\\ & \le & \frac{1}{n} \, \Big\| \alpha d_\varphi \Big(
\sum_{k=1}^n \mathsf{E}_{\mathcal{N}} \big( a_k a_k^* \big) \Big)
d_{\varphi} \alpha^* \Big\|_{L_{1/2}(\mathcal{N})}^{1/2}.
\end{eqnarray*}
This gives $$\big\| \omega_r(x) \big\|_{R_1^n(L_1(\mathcal{A}))}
\le \frac{1}{\sqrt{n}} \, \|x\|_{L_1^r(\mathcal{M},
\mathsf{E})}.$$ On the other hand, by the triangle inequality
$$\big\| \omega_r(x) \big\|_{R_1^n(L_1(\mathcal{A}))} = \frac{1}{n}
\, \Big\| \sum_{k=1}^n x_k \otimes e_{1k}
\Big\|_{R_1^n(L_1(\mathcal{A}))} \le \|x\|_{L_1(\mathcal{M})}.$$
These estimates show that $\omega_r$ is a contraction. Note also
that $$\big\langle u_r(x), \omega_r(y) \big\rangle = \frac{1}{n}
\sum_{k=1}^n \mbox{tr}_{\mathcal{A}} \big( x_{k}^* y_{k}^{} \big)
= \frac{1}{n} \sum_{k=1}^n \mbox{tr}_{\mathcal{M}} (x^* y) =
\langle x, y \rangle.$$ In particular, since it follows from
Corollary 2.12 of \cite{J1} that
$$\mathcal{R}_{\infty,1}^n (\mathcal{M}; \mathsf{E}) = \Big(
L_1(\mathcal{M}) + \frac{1}{\sqrt{n}} \, L_1^r(\mathcal{M},
\mathsf{E}) \Big)^*,$$ it turns out that the map $\omega_r^* u_r$
is the identity on $\mathcal{R}_{\infty,1}^n (\mathcal{M};
\mathsf{E})$ and $u_r \omega_r^*$ is a bounded projection onto the
image of $u_r$ with constants independent of $n$. This completes
the proof in the row case. The column case follows in the same
way. \fin

The careful reader will have observed that the projection maps
$u_r w_r^*$ and $u_c w_c^*$ are algebraically the same map $u
w^*$, modulo the identification of $R_n(\mathcal{A})$ and
$C_n(\mathcal{A})$ with $\mathcal{A}_{\oplus n}$. If we consider
the conditional expectation $\mathcal{E}_{\mathsf{A}_k}:
\mathcal{A} \to \mathsf{A}_k$ and we set
$\mathcal{E}_{\mathsf{A}_k}(z_k) = (\alpha_k,\beta_k)$, then
$$uw^* \Big( \sum_{k=1}^n z_k \otimes \delta_k \Big) = \sum_{k=1}^n
\pi_k \Big( \frac{1}{2n} \sum_{j=1}^n \alpha_j - \beta_j,
\frac{1}{2n} \sum_{j=1}^n \beta_j - \alpha_j \Big) \otimes
\delta_k.$$ In particular, Lemma
\ref{Lemma-Complemented-Isomorphism} allows us to identify the
interpolation space $$\mathrm{X}_{\frac12} = \big[
\mathcal{C}_{\infty,1}^n (\mathcal{M}, \mathsf{E}),
\mathcal{R}_{\infty,1}^n (\mathcal{M}, \mathsf{E})
\big]_{\frac12}$$ with a complemented subspace of
$L_\infty(\mathcal{A}; \mathrm{OH}_n)$. This will be implicitly
used below. However, the difficult part in proving Theorem
\ref{MAMSth} is to identify the norm of the space
$\mathrm{X}_{\frac12}$. Of course, according to the fact that we
are interpolating $2$-term intersection spaces, we expect a
$4$-term maximum. This is the case and we define
$\mathcal{J}_{\infty,2}^n(\mathcal{M}, \mathsf{E})$ as the space
of elements $x$ in $\mathcal{M}$ equipped with the norm
$$\max_{u, v \in \{4,\infty\}} \left\{ n^{\frac{1}{\xi(u,v)}} \,
\sup \big\{ \|\alpha x \beta\|_{L_{\xi(u,v)}(\mathcal{M})} \, | \,
\|\alpha\|_{L_u(\mathcal{N})}, \|\beta\|_{L_v(\mathcal{N})} \le 1
\big\} \right\},$$ where $\xi(u,v)$ is given by
$\frac{1}{\xi(u,v)} = \frac1u + \frac1v$. Obviously, multiplying
by elements $\alpha, \beta$ in the unit ball of
$L_\infty(\mathcal{N})$ and taking suprema does not contribute to
the corresponding $L_{\xi(u,v)}(\mathcal{M})$ term. In other
words, we may rewrite the norm of $x$ in $\mathcal{J}_{\infty,
2}^n(\mathcal{M}, \mathsf{E})$ as $$\|x\|_{\mathcal{J}_{\infty,
2}^n(\mathcal{M}, \mathsf{E})} = \max \Big\{
\|x\|_{\Lambda_{(u,v)}^n} \, \big| \ u,v \in \{4,\infty\} \Big\}$$
where the $\Lambda_{(u,v)}^n$ norms are given by
$$\begin{array}{lcl} \|x\|_{\Lambda_{(\infty, \infty)}^n} & = &
\displaystyle \|x\|_{\mathcal{M}}, \\ [7pt]
\|x\|_{\Lambda_{(\infty, 4)}^n} & = & n^{\frac14} \, \displaystyle
\sup \Big\{ \|x \beta\|_{L_4(\mathcal{M})} \, \big|
\ \|\beta\|_{L_4(\mathcal{N})} \le 1 \Big\}, \\
[7pt] \|x\|_{\Lambda_{(4,\infty)}^n} & = & n^{\frac14} \,
\displaystyle \sup \Big\{ \|\alpha x\|_{L_4(\mathcal{M})} \, \big|
\ \|\alpha\|_{L_4(\mathcal{N})} \le 1 \Big\}, \\
[7pt]\|x\|_{\Lambda_{(4,4)}^n} & = & n^{\frac12} \, \displaystyle
\sup \Big\{ \|\alpha x \beta\|_{L_2(\mathcal{M})} \, \big| \
\|\alpha\|_{L_4(\mathcal{N})}, \|\beta\|_{L_4(\mathcal{N})} \le 1
\Big\}. \end{array}$$ These norms arise as particular cases of the
so-called conditional $L_p$ spaces, which were defined in
\cite{JP2} and will be further exploited in \cite{JP4}. Before
identifying the norm of $\mathrm{X}_{\frac12}$, we need some
information on interpolation spaces.

\begin{lemma} \label{Lemintinfty}
If $(1/u,1/v) = (\theta/2, (1-\theta)/2)$, we have for $x \in
\mathcal{M}$
\begin{eqnarray*}
\|x\|_{[\mathcal{M}, L_{\infty}^r(\mathcal{M},
\mathsf{E})]_{\theta}} & = & \sup \Big\{ \|\alpha x
\|_{L_u(\mathcal{M})} \, \big| \
\|\alpha\|_{L_u(\mathcal{N})} \le 1 \Big\}, \\
\|x\|_{[L_{\infty}^c(\mathcal{M}, \mathsf{E}),
\mathcal{M}]_{\theta}} & = & \sup \Big\{ \|x
\beta\|_{L_v(\mathcal{M})} \, \big| \
\|\beta\|_{L_v(\mathcal{N})} \le 1 \Big\}, \\
\|x\|_{[L_{\infty}^c(\mathcal{M}, \mathsf{E}),
L_{\infty}^r(\mathcal{M}, \mathsf{E})]_{\theta}} & = & \sup \Big\{
\|\alpha x \beta\|_{L_2(\mathcal{M})} \, \big| \
\|\alpha\|_{L_u(\mathcal{N})}, \|\beta\|_{L_v(\mathcal{N})} \le 1
\Big\}.
\end{eqnarray*}
\end{lemma}

The proof can be found in \cite{JP2}. In the finite setting, this
result follows from a well-known application of
Helson/Lowdenslager, Wiener/Masani type results on the existence
of operator-valued analytic functions. This kind of applications
has been used extensively by Pisier in his theory of vector-valued
$L_p$ spaces. We shall also need an explicit expression for the
norm of $L_\infty(\mathcal{A}; \mathrm{OH}_n)$. The formula below
was originally obtained by Pisier in \cite{P0}. The proof for
non-hyperfinite von Neumann algebras (as it is our case) can be
found in \cite{JP2,X}.

\begin{lemma} \label{pisint}
We have $$\Big\| \sum_{k=1}^n a_k \ten \delta_k
\Big\|_{L_\infty(\mathcal{A}; \mathrm{OH}_n)} = \sup \left\{
\Big\| \sum_{k=1}^n a_k^* \alpha a_k
\Big\|_{L_2(\mathcal{A})}^{\frac12} \, \big| \ \alpha \ge 0, \,
\|\alpha\|_2 \le 1 \right\}.$$
\end{lemma}

\noindent Now we are ready for the main result in this section.

\begin{theorem} \label{Theorem-Intersection} We have
isomorphically
\begin{eqnarray*}
\big[ \mathcal{C}_{\infty,1}^n (\mathcal{M}, \mathsf{E}),
\mathcal{R}_{\infty,1}^n (\mathcal{M}, \mathsf{E}) \big]_{\frac12}
& \simeq & \mathcal{J}_{\infty,2}^n (\mathcal{M}, \mathsf{E}).
\end{eqnarray*}
Moreover, the constants in these isomorphisms are uniformly
bounded on $n$.
\end{theorem}

\dem We have a contractive inclusion
$$\mathrm{X}_{\frac12} \subset \mathcal{J}_{\infty,2}^n
(\mathcal{M}, \mathsf{E}).$$ Indeed, by elementary properties of
interpolation spaces we find $$\mathrm{X}_{\frac12} \subset
[\mathcal{M}, \mathcal{M}]_{\frac12} \cap [\sqrt{n} L_\infty^c,
\mathcal{M}]_{\frac12} \cap [\mathcal{M}, \sqrt{n}
L_\infty^r]_{\frac12} \cap [\sqrt{n} L_\infty^c, \sqrt{n}
L_\infty^r]_{\frac12},$$ where $L_\infty^r$ and $L_\infty^c$ are
abbreviations for $L_\infty^r(\mathcal{M}, \mathsf{E})$ and
$L_\infty^c(\mathcal{M}, \mathsf{E})$ respectively. Using the
obvious identity $[\lambda_0 \mathrm{X}_0, \lambda_1
\mathrm{X}_1]_\theta = \lambda_0^{1-\theta} \lambda_1^\theta \,
\mathrm{X}_\theta$ and applying Lemma \ref{Lemintinfty} we
rediscover the norm of the space $\mathcal{J}_{\infty,2}^n
(\mathcal{M}, \mathsf{E})$ on the right hand side. This gives that
the inclusion $\mathrm{X}_{\frac12} \subset
\mathcal{J}_{\infty,2}^n (\mathcal{M}, \mathsf{E})$ holds
contractively.

\vskip5pt

\noindent To prove the reverse inclusion, we note from Lemmas
\ref{Lemma-Complemented-Isomorphism} and \ref{pisint} that
\begin{eqnarray*}
\|x\|_{\mathrm{X}_{\frac12}} & \sim & \big\| u(x)
\big\|_{[C_n(\mathcal{A}), R_n(\mathcal{A})]_{\frac{1}{2}}} \\ & =
& \Big\| \sum_{k=1}^n x_k \ten \delta_k
\Big\|_{L_\infty(\mathcal{A}; \mathrm{OH}_n)} \\ & = & \sup
\left\{ \Big\| \sum_{k=1}^n x_k^* a x_k
\Big\|_{L_2(\mathcal{A})}^{\frac12} \, \big| \ a \ge 0, \, \|a\|_2
\le 1 \right\} = \mathrm{A}.
\end{eqnarray*}
Thus, it remains to see that
\begin{equation} \label{Equation-Remaining-Estimate2}
\mathrm{A} \ \lesssim \max \Big\{ \|x\|_{\Lambda_{(u,v)}^n} \,
\big| \ u,v \in \{4,\infty\} \Big\} =
\|x\|_{\mathcal{J}_{\infty,2}^n(\mathcal{M}, \mathsf{E})}.
\end{equation}
In order to justify \eqref{Equation-Remaining-Estimate2}, we
introduce the orthogonal projections $\mathsf{L}_k$ and
$\mathsf{R}_k$ in $L_2(\mathcal{A})$ defined as follows. Given $a
\in L_2(\mathcal{A})$, the vector $\mathsf{L}_k(a)$ (resp.
$\mathsf{R}_k(a)$) collects the reduced words in $a$ starting
(resp. ending) with a letter in $\mathsf{A}_k$. In other words,
following standard terminology in free probability, we have
\begin{eqnarray*}
\mathsf{L}_k: L_2(\mathcal{A}) & \longrightarrow & L_2 \Big( \big[
\bigoplus_{m \ge 1} \bigoplus_{j_1 =k \neq j_2 \neq \cdots
\neq j_m} \bubl_{j_1} \bubl_{j_2} \cdots \bubl_{j_m} \big]'' \Big), \\
\mathsf{R}_k: L_2(\mathcal{A}) & \longrightarrow & L_2 \Big( \big[
\bigoplus_{m \ge 1} \bigoplus_{j_1 \neq j_2 \neq \cdots \neq k =
j_m} \bubl_{j_1} \bubl_{j_2} \cdots \bubl_{j_m} \big]'' \Big).
\end{eqnarray*}
Now, given a positive operator $a$ in $L_2(\mathcal{A})$ and a
fixed integer $1 \le k \le n$, we consider the following way to
decompose $a$ in terms of the projections $\mathsf{L}_k$ and
$\mathsf{R}_k$ and the conditional expectation
$\mathsf{E}_\mathcal{N}: \mathcal{A} \to \mathcal{N}$
\begin{equation} \label{Equation-Decomposition}
a = \mathsf{E}_{\mathcal{N}}(a) + \big( \mathsf{L}_k(a) +
\mathsf{R}_k(a) - \mathsf{R}_k \mathsf{L}_k(a) \big) +
\gamma_k(a),
\end{equation} where the term $\gamma_k(a)$ has the
following form $$\gamma_k(a) = a - \mathsf{E}_{\mathcal{N}}(a) -
\mathsf{L}_k(a) - \mathsf{R}_k(a - \mathsf{L}_k(a)).$$ The
triangle inequality gives $\mathrm{A}^2 \le \mathrm{A}_1^2 +
\mathrm{A}_2^2 + \mathrm{A}_3^2$, where the terms $\mathrm{A}_j$
are the result of replacing $a$ in $\mathrm{A}$ by the $j$-th term
in the decomposition \eqref{Equation-Decomposition}. The terms in
the bracket of \eqref{Equation-Decomposition} are understood as
one term which gives rise to $\mathrm{A}_2$. For the first term we
use
\begin{eqnarray*}
\mathrm{A}_1^2 & = & \Big\| \sum_{k=1}^n x_k^*
\mathsf{E}_{\mathcal{N}}(a) x_k \Big\|_{L_2(\mathcal{A})} \\
& \le & \Big\| \sum_{k=1}^n \mathsf{E}_{\mathcal{N}} \big( x_k^*
\mathsf{E}_{\mathcal{N}}(a) x_k \big) \Big\|_{L_2(\mathcal{N})} \\
& + & \Big\| \sum_{k=1}^n x_k^* \mathsf{E}_{\mathcal{N}}(a) x_k -
\mathsf{E}_{\mathcal{N}} \big( x_k^* \mathsf{E}_{\mathcal{N}}(a)
x_k \big) \Big\|_{L_2(\mathcal{A})} = \mathrm{A}_{11}^2 +
\mathrm{A}_{12}^2.
\end{eqnarray*}
Since $\mathsf{E}_{\mathcal{N}} \big( x_k^*
\mathsf{E}_{\mathcal{N}}(a) x_k \big) = \mathsf{E} \big( x^*
\mathsf{E}_{\mathcal{N}}(a) x \big)$ and $a \in
\mathsf{B}_{L_2(\mathcal{A})}^+$, we obtain
\begin{eqnarray*}
\mathrm{A}_{11} & = & n^{\frac12} \, \sup \Big\{
\mbox{tr}_{\mathcal{N}} \Big( \beta^* \mathsf{E} \big( x^*
\mathsf{E}_{\mathcal{N}}(a) x \big) \beta \Big)^{\frac12} \, \big|
\ \|\beta\|_{L_4(\mathcal{N})} \le 1 \Big\} \\ & \le & n^{\frac12}
\, \sup \Big\{ \mbox{tr}_{\mathcal{M}} \big( \beta^* x^* \alpha^*
\alpha x \beta \big)^{\frac12} \, \big| \
\|\alpha\|_{L_4(\mathcal{N})}, \|\beta\|_{L_4(\mathcal{N})} \le 1
\Big\}.
\end{eqnarray*}
This gives $\mathrm{A}_{11} \le \|x\|_{\Lambda_{(4,4)}^n} \le
\|x\|_{\mathcal{J}_{\infty,2}^n(\mathcal{M}, \mathsf{E})}$. On the
other hand, by freeness
\begin{eqnarray*}
\mathrm{A}_{12}^2 & = & \Big( \sum_{k=1}^n \big\| x_k^*
\mathsf{E}_{\mathcal{N}}(a) x_k - \mathsf{E}_{\mathcal{N}} \big(
x_k^* \mathsf{E}_{\mathcal{N}}(a) x_k \big)
\big\|_{L_2(\mathcal{A})}^2 \Big)^{\frac12} \\ & \le & 2 \, \Big(
\sum_{k=1}^n \big\| x_k^* \mathsf{E}_{\mathcal{N}}(a) x_k
\big\|_{L_2(\mathcal{A})}^2 \Big)^{\frac12} = 2 \, n^{\frac12} \,
\big\| x^* \mathsf{E}_{\mathcal{N}}(a) x
\big\|_{L_2(\mathcal{M})}.
\end{eqnarray*}
Then positivity gives
$$\mathrm{A}_{12} \le \sqrt{2} \, \|x\|_{\Lambda_{(4,\infty)}^n}
\le \sqrt{2} \, \|x\|_{\mathcal{J}_{\infty,2}^n(\mathcal{M},
\mathsf{E})}.$$ Let us now estimate the term $\mathrm{A}_2$
\begin{eqnarray*}
\lefteqn{\Big\| \sum_{k=1}^n x_k^* \big( \mathsf{L}_k(a) +
\mathsf{R}_k(a) - \mathsf{R}_k \mathsf{L}_k(a) \big) x_k
\Big\|_{L_2(\mathcal{A})}}
\\ & = & \sup \left\{ \sum_{k=1}^n \mbox{tr}_{\mathcal{A}} \Big( b
x_k^* \big( \mathsf{L}_k(a) + \mathsf{R}_k(a) - \mathsf{R}_k
\mathsf{L}_k(a) \big) x_k \Big) \, \big| \
\|b\|_{L_2(\mathcal{A})} \le 1 \right\} \\ & = & \sup \left\{
\sum_{k=1}^n \mbox{tr}_{\mathcal{A}} \Big( x_k b x_k^* \big(
\mathsf{L}_k(a) + \mathsf{R}_k(a) - \mathsf{R}_k \mathsf{L}_k(a)
\big) \Big) \, \big| \ \|b\|_{L_2(\mathcal{A})} \le 1 \right\}
\\ & \le & \sup_{\|b\|_{L_2(\mathcal{A})} \le 1} \Big( \sum_{k=1}^n
\big\| x_k b x_k^* \big\|_{L_2(\mathcal{A})}^2 \Big)^{\frac12}
\Big( \sum_{k=1}^n \big\| \mathsf{L}_k(a) + \mathsf{R}_k(a) -
\mathsf{R}_k \mathsf{L}_k(a) \big\|_{L_2(\mathcal{A})}^2
\Big)^{\frac12}.
\end{eqnarray*}
The second factor on the right is estimated by orthogonality
$$\Big( \sum_{k=1}^n \big\| \mathsf{L}_k(a) + \mathsf{R}_k(a) -
\mathsf{R}_k \mathsf{L}_k(a) \big\|_{L_2(\mathcal{A})}^2
\Big)^{\frac12} \le 3 \, \|a\|_{L_2(\mathcal{A})} \le 3.$$ For the
first factor, we write $b$ as a linear combination $(b_1 - b_2) +
i (b_3 - b_4)$ of four positive operators. Therefore, all these
terms are covered by the following estimate, to be proved below.

\vskip5pt

\noindent \textbf{Claim.} Given $a \in
\mathsf{B}_{L_2(\mathcal{A})}^+$, we have
\begin{equation} \label{claim}
\Big( \sum_{k=1}^n \big\| x_k a x_k^* \big\|_{L_2(\mathcal{A})}^2
\Big)^{\frac14} \lesssim \max \Big\{ \|x\|_\mathcal{M}, \,
n^{\frac14} \sup_{\|\beta\|_{L_4(\mathcal{N})} \le 1} \|x
\beta\|_{L_4(\mathcal{M})} \Big\}.
\end{equation}
Before justifying our claim, we complete the proof. It remains to
estimate the term $\mathrm{A}_3$ associated to $\gamma_k(a)$. We
first observe that $\gamma_k(a)$ is a mean-zero element of
$L_2(\mathcal{A})$ made up of reduced words not starting nor
ending with a letter in $\mathsf{A}_k$. Indeed, note that
$\mathsf{E}_{\mathcal{N}}(\gamma_k(a)) = 0$ and that we first
eliminate the words starting with a letter in $\mathsf{A}_k$ by
subtracting $\mathsf{L}_k(a)$ and, after that, we eliminate the
remaining words which end with a letter in $\mathsf{A}_k$ by
subtracting $\mathsf{R}_k(a - \mathsf{L}_k(a))$. Therefore, it
turns out that the family of random variables $x_1^* \gamma_1(a)
x_1, x_2^* \gamma_2(a) x_2, \ldots, x_n^* \gamma_n(a) x_n$ is free
over $\mathcal{N}$. In particular, by orthogonality
\begin{eqnarray*}
\Big\| \sum_{k=1}^n x_k^* \gamma_k(a) x_k
\Big\|_{L_2(\mathcal{A})}^{\frac12} & = & \Big( \sum_{k=1}^n
\big\| x_k^* \gamma_k(a) x_k \big\|_{L_2(\mathcal{A})}^2
\Big)^{\frac14}.
\end{eqnarray*}
However, recalling that $\gamma_k(a)$ is a mean-zero element made
up of words not starting nor ending with a letter in
$\mathsf{A}_k$, the following identities hold for the conditional
expectation $\mathcal{E}_{\mathsf{A}_k}: L_1(\mathcal{A}) \to
L_1(\mathsf{A}_k)$
\begin{equation} \label{Expectation-Only}
\begin{array}{c}
\mathcal{E}_{\mathsf{A}_k} \Big( \gamma_k(a)^* \big( x_kx_k^* -
\mathsf{E}_{\mathcal{N}}(x_kx_k^*) \big) \gamma_k(a) \Big) = 0, \\
[5pt] \mathcal{E}_{\mathsf{A}_k} \Big( \gamma_k(a) \big( x_kx_k^*
- \mathsf{E}_{\mathcal{N}}(x_kx_k^*) \big) \gamma_k(a)^* \Big) =
0.
\end{array}
\end{equation}
Using property (\ref{Expectation-Only}) we find
\begin{eqnarray*}
\big\|x_k^* \gamma_k(a) x_k \big\|_{L_2(\mathcal{A})}^2 & = &
\mbox{tr}_{\mathcal{A}} \big( x_k^* \gamma_k(a)^* x_k x_k^*
\gamma_k(a) x_k \big) \\ & = & \mbox{tr}_{\mathcal{A}} \big( x_k^*
\mathcal{E}_{\mathsf{A}_k} \big( \gamma_k(a)^* x_k x_k^*
\gamma_k(a) \big) x_k \big) \\ & = & \mbox{tr}_{\mathcal{A}} \big(
x_k^* \mathcal{E}_{\mathsf{A}_k} \big( \gamma_k(a)^*
\mathsf{E}_{\mathcal{N}}(x_k x_k^*) \gamma_k(a) \big) x_k \big) \\
& = & \mbox{tr}_{\mathcal{A}} \big( \gamma_k(a) x_k x_k^*
\gamma_k(a)^* \mathsf{E}_{\mathcal{N}}(x_k x_k^*) \big) \\ & = &
\mbox{tr}_{\mathcal{A}} \big( \mathcal{E}_{\mathsf{A}_k} \big(
\gamma_k(a) x_k x_k^* \gamma_k(a)^* \big)
\mathsf{E}_{\mathcal{N}}(x_k x_k^*) \big) \\ & = &
\mbox{tr}_{\mathcal{A}} \big( \mathcal{E}_{\mathsf{A}_k} \big(
\gamma_k(a) \mathsf{E}_{\mathcal{N}}(x_k x_k^*) \gamma_k(a)^*
\big) \mathsf{E}_{\mathcal{N}}(x_k x_k^*) \big) \\ & = &
\mbox{tr}_{\mathcal{A}} \big( \mathsf{E}_{\mathcal{N}}(x_k
x_k^*)^{\frac12} \gamma_k(a) \mathsf{E}_{\mathcal{N}}(x_k x_k^*)
\gamma_k(a)^* \mathsf{E}_{\mathcal{N}}(x_k x_k^*)^{\frac12} \big).
\end{eqnarray*}
In conjunction with $\|\gamma_k(a)\|_2 \le 5 \, \|a\|_2$ and
H{\"o}lder's inequality, this yields
$$\big\|x_k^* \gamma_k(a) x_k \big\|_{L_2(\mathcal{A})}^2 = \big\|
\mathsf{E}_{\mathcal{N}}(x x^*)^{\frac12} \gamma_k(a)
\mathsf{E}_{\mathcal{N}}(x x^*)^{\frac12}
\big\|_{L_2(\mathcal{A})}^2 \le 25 \,
\|x\|_{L_{\infty}^r(\mathcal{M},\mathsf{E})}^4.$$ We refer to
Lemma \ref{Lemintinfty} or \cite{JP2} for the fact that
$$\|x\|_{L_\infty^r(\mathcal{M}, \mathsf{E})} \le \sup \Big\{
\|\alpha x\|_{L_4(\mathcal{M})} \, \big| \
\|\alpha\|_{L_4(\mathcal{N})} \le 1 \Big\}.$$ The inequalities
proved so far give rise to the following estimate $$\Big\|
\sum_{k=1}^n x_k^* \gamma_k(a) x_k
\Big\|_{L_2(\mathcal{A})}^{\frac12} \le \sqrt{5} \,
\|x\|_{\Lambda_{(4,\infty)}^n} \le \sqrt{5} \,
\|x\|_{\mathcal{J}_{\infty,2}^n(\mathcal{M}, \mathsf{E})}.$$
Therefore, it remains to prove the claim. We proceed in a similar
way. According to the decomposition
\eqref{Equation-Decomposition}, we may use the triangle inequality
and decompose the left hand side of \eqref{claim} into three terms
$\mathrm{B}_1, \mathrm{B}_2, \mathrm{B}_3$. For the first term, we
deduce from positivity that $$\Big( \sum_{k=1}^n \big\| x_k
\mathsf{E}_\mathcal{N}(a) x_k^* \big\|_{L_2(\mathcal{A})}^2
\Big)^{\frac14} = n^{\frac14} \big\| x \mathsf{E}_\mathcal{N}(a)
x^* \big\|_{L_2(\mathcal{M})}^{\frac12} \le n^{\frac14}
\sup_{\|\beta\|_{L_4(\mathcal{N})} \le 1} \|x
\beta\|_{L_4(\mathcal{M})}.$$ On the other hand, it is elementary
that
\begin{eqnarray*}
\lefteqn{\hskip-15pt \Big( \sum_{k=1}^n \big\| x_k^* \big(
\mathsf{L}_k(a) + \mathsf{R}_k(a) - \mathsf{R}_k \mathsf{L}_k(a)
\big) x_k \big\|_{L_2(\mathcal{A})}^2 \Big)^{\frac14}} \\ & \le &
\|x\|_\mathcal{M} \Big( \sum_{k=1}^n \big\| x_k^* \big(
\mathsf{L}_k(a) + \mathsf{R}_k(a) - \mathsf{R}_k \mathsf{L}_k(a)
\big) x_k \big\|_{L_2(\mathcal{A})}^2 \Big)^{\frac14} \le 3 \,
\|x\|_\mathcal{M}.
\end{eqnarray*}
This leaves us with the term $\mathrm{B}_3$. Arguing as above
$$\Big( \sum_{k=1}^n \big\| x_k \gamma_k(a) x_k^*
\big\|_{L_2(\mathcal{A})}^2 \Big)^{\frac14} \le \sqrt{5} \
n^{\frac14} \sup_{\|\beta\|_{L_4(\mathcal{N})} \le 1} \|x
\beta\|_{L_4(\mathcal{M})}.$$ Therefore, the claim holds and the
proof is complete. \fin

\begin{remark}
\emph{The arguments in Theorem \ref{Theorem-Intersection} also
give
\begin{eqnarray*}
\|x\|_{[\mathcal{M}, \mathcal{R}_{\infty,1}^n (\mathcal{M},
\mathsf{E})]_{\frac12}} & \sim & \max \Big\{ \|x\|_\mathcal{M},
\|x\|_{\Lambda_{(4,\infty)}^n} \Big\}, \\
\|x\|_{[\mathcal{C}_{\infty,1}^n (\mathcal{M}, \mathsf{E}) \, , \,
\mathcal{M}]_{\frac12}} & \sim & \max \Big\{ \|x\|_{\mathcal{M}},
\|x\|_{\Lambda_{(\infty,4)}^n} \Big\}.
\end{eqnarray*}}
\end{remark}

Now we show how the space $\mathrm{X}_{1/2}$ is related to Theorem
\ref{MAMSth}. The idea follows from a well-known argument in which
complete boundedness arises as a particular case of amalgamation.
More precisely, if $L_2^r(\mathcal{M})$/$L_2^c(\mathcal{M})$
denote the row/column quantizations of $L_2(\mathcal{M})$ and $2
\le q \le \infty$, the row/column operator space structures on
$L_q(\mathcal{M})$ are defined as follows
\begin{equation} \label{Eq-RCLp}
\begin{array}{rcl}
L_q^r(\mathcal{M}) & = & \big[ \mathcal{M},
L_2^r(\mathcal{M})\big]_{\frac{2}{q}}, \\ [5pt] L_q^c(\mathcal{M})
& = & \big[ \mathcal{M}, L_2^c(\mathcal{M})\big]_{\frac{2}{q}}.
\end{array}
\end{equation}
The following result from \cite{JP2} is a generalized form of
\eqref{osstruct} in the Introduction.

\begin{lemma} \label{Lemma-Isometry-Rho}
If $\mathcal{M}_m = \mathrm{M}_m(\mathcal{M})$, we have
\begin{eqnarray*}
\big\| d_{\varphi}^{\frac{1}{4}} \big( x_{ij} \big)
\big\|_{\mathrm{M}_m(L_4^r(\mathcal{M}))} & = &
\sup_{\|\alpha\|_{S_4^m} \le 1} \Big\| d_\varphi^{\frac{1}{4}}
\Big( \sum_{k=1}^m \alpha_{ik} x_{kj} \Big)
\Big\|_{L_4(\mathcal{M}_m)}, \\ \big\| \big( x_{ij} \big)
d_{\varphi}^{\frac{1}{4}}
\big\|_{\mathrm{M}_m(L_4^c(\mathcal{M}))} & = &
\sup_{\|\beta\|_{S_4^m} \le 1} \Big\| \Big( \sum_{k=1}^m x_{ik}
\beta_{kj} \Big) d_\varphi^{\frac{1}{4}}
\Big\|_{L_4(\mathcal{M}_m)}.
\end{eqnarray*}
\end{lemma}

\noindent The proof follows from
\begin{equation} \label{Eq-Schatten-RC}
\begin{array}{rcl}
\big\| d_{\varphi}^{\frac{1}{2}} \big( x_{ij} \big)
\big\|_{\mathrm{M}_m(L_2^r(\mathcal{M}))} & = & \displaystyle
\sup_{\|\alpha\|_{S_2^m} \le 1} \Big\| d_\varphi^{\frac{1}{2}}
\Big( \sum_{k=1}^m \alpha_{ik} x_{kj} \Big)
\Big\|_{L_2(\mathcal{M}_m)}, \\ \big\| \big( x_{ij} \big)
d_{\varphi}^{\frac{1}{2}}
\big\|_{\mathrm{M}_m(L_2^c(\mathcal{M}))} & = & \displaystyle
\sup_{\|\beta\|_{S_2^m} \le 1} \Big\| \Big( \sum_{k=1}^m x_{ik}
\beta_{kj} \Big) d_\varphi^{\frac{1}{2}}
\Big\|_{L_2(\mathcal{M}_m)},
\end{array}
\end{equation}
and some complex interpolation formulas developed in \cite{JP2}.
The identity \eqref{Eq-Schatten-RC} from which we interpolate is a
well-known expression in operator space theory, see e.g. p.56 in
\cite{ER}. Now we define the space
$\mathcal{J}_{\infty,2}^n(\mathcal{M})$ as follows
$$\mathcal{J}_{\infty,2}^n(\mathcal{M}) = \mathcal{M} \cap
n^{\frac{1}{4}} L_4^c(\mathcal{M}) \cap n^{\frac{1}{4}}
L_4^r(\mathcal{M}) \cap n^{\frac{1}{2}} L_2(\mathcal{M}).$$ Lemma
\ref{Lemma-Isometry-Rho} determines the operator space structure
of the cross terms in $\mathcal{J}_{\infty,2}^n(\mathcal{M})$. On
the other hand, according to Pisier's fundamental identity
\eqref{fbasexu} or to Lemma \ref{pisint}, it is easily seen that
we have
$$\big\| d_{\varphi}^{\frac{1}{4}} \big( x_{ij} \big)
d_{\varphi}^{\frac{1}{4}} \big\|_{\mathrm{M}_m(L_2(\mathcal{M}))}
= \sup_{\|\alpha\|_{S_4^m}, \|\beta\|_{S_4^m} \le 1} \Big\|
d_\varphi^{\frac{1}{4}} \Big( \sum_{k,l=1}^m \alpha_{ik} x_{kl}
\beta_{lj} \Big) d_\varphi^{\frac{1}{4}}
\Big\|_{L_2(\mathcal{M}_m)}.$$ In other words, the o.s.s. of
$\mathcal{J}_{\infty,2}^n(\mathcal{M})$ is described by the
isometry
\begin{equation} \label{ossJ}
\mathrm{M}_m \big( \mathcal{J}_{\infty,2}^n(\mathcal{M}) \big) =
\mathcal{J}_{\infty,2}^n(\mathcal{M}_m, \mathsf{E}_m),
\end{equation}
where $\mathcal{M}_m = \mathrm{M}_m(\mathcal{M})$ and
$\mathsf{E}_m = id_{\mathrm{M}_m} \otimes \varphi: \mathcal{M}_m
\to \mathrm{M}_m$ for $m \ge 1$. This means that the
\emph{vector-valued} spaces $\mathcal{J}_{\infty,2}^n(\mathcal{M},
\mathsf{E})$ describe the o.s.s. of the \emph{scalar-valued}
spaces $\mathcal{J}_{\infty,2}^n(\mathcal{M})$. In the result
below we prove the operator space/free analogue of a form of
Rosenthal's inequality in the limit case $p \to \infty$, see
\cite{JP2,JP4} for more details. This result does not have a
commutative counterpart. The particular case for $\mathcal{M} =
\mathcal{B}(\ell_2^n)$ recovers Theorem \ref{MAMSth}. Given a von
Neumann algebra $\mathcal{M}$, we set as above $\mathsf{A}_k =
\mathcal{M} \oplus \mathcal{M}$.

\begin{corollary} \label{Corollary-Sigma-Infty-2}
If $\mathcal{A}_\mathcal{N} = *_{\mathcal{N}} \mathsf{A}_k$, the
map $$u: x \in \mathcal{J}_{\infty,2}^n(\mathcal{M}, \mathsf{E})
\mapsto \sum_{k=1}^n x_k \otimes \delta_k \in
L_\infty(\mathcal{A}_\mathcal{N}; \mathrm{OH}_n)$$ is an
isomorphism with complemented image and constants independent of
$n$. In particular, replacing as usual $(\mathcal{M}, \mathcal{N},
\mathsf{E})$ by $(\mathcal{M}_m, \mathrm{M}_m, \mathsf{E}_m)$ and
replacing $\mathcal{A}_\mathcal{N}$ by the non-amalgamated algebra
$\mathcal{A} = \mathsf{A}_1
* \mathsf{A}_2 * \cdots * \mathsf{A}_n$, we obtain a
cb-isomorphism with cb-complemented image and constants
independent of $n$ $$\sigma: x \in
\mathcal{J}_{\infty,2}^n(\mathcal{M}) \mapsto \sum_{k=1}^n x_k
\otimes \delta_k \in L_\infty(\mathcal{A}; \mathrm{OH}_n).$$
\end{corollary}

\dem The first assertion follows from Lemma
\ref{Lemma-Complemented-Isomorphism} and Theorem
\ref{Theorem-Intersection}. To prove the second assertion we
choose the triple $(\mathcal{M}_m, \mathrm{M}_m, \mathsf{E}_m)$
and apply \eqref{ossJ}. This provides us with an isomorphic
embedding $$\sigma_m: x \in \mathrm{M}_m \big(
\mathcal{J}_{\infty,2}^n (\mathcal{M}) \big) \mapsto \sum_{k=1}^n
x_k \otimes \delta_k \in L_\infty (\mathcal{A}_m;
\mathrm{OH}_n),$$ where the von Neumann algebra $\mathcal{A}_m$ is
given by $$\mathcal{A}_m = \mathrm{M}_m(\mathcal{A}) =
\mathrm{M}_m(\mathsf{A}_1) *_{\mathrm{M}_m}
\mathrm{M}_m(\mathsf{A}_2) *_{\mathrm{M}_m} \cdots
*_{\mathrm{M}_m} \mathrm{M}_m(\mathsf{A}_n).$$ The last isometry
is well-known, see e.g. \cite{J2}. In particular
$$L_\infty(\mathcal{A}_m; \mathrm{OH}_n) = \mathrm{M}_m \big(
L_\infty(\mathcal{A}; \mathrm{OH}_n) \big)$$ and it turns out that
the map $\sigma_m = id_{\mathrm{M}_m} \otimes \sigma$. This
completes the proof. \fin

\section{Complete embedding of Schatten classes}
\label{S3CBE}

Given $1 < q \le 2$, we construct a completely isomorphic
embedding of $S_q$ into the predual of a $\mathrm{QWEP}$ von
Neumann algebra. In fact, we prove Theorem \ref{QSQSQS} and deduce
our cb-embedding via Lemma \ref{Lemma-Motivation}. Theorem
\ref{lp} and the subsequent family of operator space $q$-stable
random variables arise by injecting the space $\ell_q$ into the
diagonal of the Schatten class $S_q$.

\subsection{Embedding of $\mathcal{K}$-spaces}

We fix $\mathcal{M} = \mathcal{B}(\ell_2)$ and consider a family
$\gamma_1, \gamma_2, \ldots$ of strictly positive numbers. Then we
consider the diagonal operator on $\ell_2$ defined by
$\mathsf{d}_\gamma = \sum_k \gamma_k e_{kk}$. This operator can be
regarded as the density $d_{\psi}$ associated to a normal strictly
semifinite faithful (\emph{n.s.s.f.} in short) weight $\psi$ on
$\mathcal{B}(\ell_2)$. Let us set $q_n$ to be the projection
$\sum_{k \le n} e_{kk}$ and let us consider the restriction of
$\psi$ to the subalgebra $q_n \mathcal{B}(\ell_2) q_n$ $$\psi_n
\Big( q_n \big( \summ_{ij} x_{ij} e_{ij} \big) q_n \Big) =
\sum_{k=1}^n \gamma_k x_{kk}.$$ Note that if we set $\mathrm{k}_n
= \psi_n(q_n)$, we obtain $\psi_n = \mathrm{k}_n \varphi_n$ for
some state $\varphi_n$ on $q_n \mathcal{B}(\ell_2) q_n$. If
$d_{\psi_n}$ denotes the density on $q_n \mathcal{B}(\ell_2) q_n$
associated to the weight $\psi_n$, we define the space
$\mathcal{J}_{\infty,2}(\psi_n)$ as the subspace
$$\Big\{ \big( z, z d_{\psi_n}^{\frac14}, d_{\psi_n}^{\frac14} z,
d_{\psi_n}^{\frac14} z d_{\psi_n}^{\frac14} \big) \, \big| \ z \in
q_n \mathcal{B}(\ell_2) q_n \Big\}$$ of the direct sum
$$\mathcal{L}_{\infty}^n = \big( C_n \otimes_h R_n \big)
\oplus_2 \big( C_n \otimes_h \mathrm{OH}_n \big) \oplus_2 \big(
\mathrm{OH}_n \otimes_h R_n \big) \oplus_2 \big( \mathrm{OH}_n
\otimes_h \mathrm{OH}_n \big).$$ In other words, we obtain the
intersection space considered in the Introduction
$$(C_n \ten_h R_n) \cap (C_n \ten_h \mathrm{OH}_n)
d_{\psi_n}^{\frac14} \cap d_{\psi_n}^{\frac14} (\mathrm{OH}_n
\ten_h R_n) \cap d_{\psi_n}^{\frac14} (\mathrm{OH}_n \ten_h
\mathrm{OH}_n) d_{\psi_n}^{\frac14}.$$

\begin{lemma} \label{Lemma-PredualMAMS11}
Let us consider $$\mathcal{K}_{1,2}(\psi_n) =
\mathcal{J}_{\infty,2}(\psi_n)^*.$$ Assume that $\mathrm{k}_n =
\sum_{k=1}^n \gamma_k$ is an integer and define $\mathcal{A}_n$ to
be the $\mathrm{k}_n$-fold reduced free product of $q_n
\mathcal{B}(\ell_2) q_n \oplus q_n \mathcal{B}(\ell_2) q_n$. If
$\pi_j: q_n \mathcal{B}(\ell_2) q_n \oplus q_n \mathcal{B}(\ell_2)
q_n \to \mathcal{A}_n$ is the natural embedding into the $j$-th
component of $\mathcal{A}_n$ and $x_j = \pi_j(x,-x)$
$$\omega: x \in \mathcal{K}_{1,2}(\psi_n) \mapsto \frac{1}{\
\mathrm{k}_n} \sum_{j=1}^{\mathrm{k}_n} x_j \otimes \delta_j \in
L_1(\mathcal{A}_n; \mathrm{OH}_{\mathrm{k}_n})$$ is a cb-embedding
with cb-complemented image and constants independent of $n$.
\end{lemma}

\dem We claim that $$\mathcal{J}_{\infty,2}(\psi_n) =
\mathcal{J}_{\infty,2}^{\mathrm{k}_n} (q_n \mathcal{B}(\ell_2)
q_n)$$ completely isometrically. Indeed, by \eqref{Eq-RCLp}
\begin{eqnarray*}
\mathrm{k}_n^{\frac14} L_4^r(q_n \mathcal{B}(\ell_2) q_n,
\varphi_n) & = & \mathrm{k}_n^{\frac14} \, \big[
\mathcal{B}(\ell_2^n), L_2^r (\mathcal{B}(\ell_2^n), \varphi_n)
\big]_{\frac12} \\ & = & \mathrm{k}_n^{\frac14} \, \big[
\mathcal{B}(\ell_2^n), d_{\varphi_n}^{\frac12} L_2^r
(\mathcal{B}(\ell_2^n), \mathrm{tr}_n) \big]_{\frac12} \\ & = &
\mathrm{k}_n^{\frac14} d_{\varphi_n}^{\frac14} \big[ C_n \ten_h
R_n, R_n \ten_h R_n \big]_{\frac12} = d_{\psi_n}^{\frac14}
(\mathrm{OH}_n \ten_h R_n).
\end{eqnarray*}
We can treat the other terms similarly and obtain
\begin{eqnarray*}
\mathrm{k}_n^{\frac14} L_4^c(q_n \mathcal{B}(\ell_2) q_n,
\varphi_n) & = & (C_n \ten_h \mathrm{OH}_n) d_{\psi_n}^{\frac14},
\\ \mathrm{k}_n^{\frac12} L_2(q_n \mathcal{B}(\ell_2) q_n,
\varphi_n) & = & d_{\psi_n}^{\frac14} (\mathrm{OH}_n \ten_h
\mathrm{OH}_n) d_{\psi_n}^{\frac14}.
\end{eqnarray*}
In particular, Corollary \ref{Corollary-Sigma-Infty-2} provides a
cb-isomorphism
$$\sigma: x \in \mathcal{J}_{\infty,2}(\psi_n) \mapsto
\sum_{j=1}^{\mathrm{k}_n} x_j \otimes \delta_j \in
L_{\infty}(\mathcal{A}_n; \mathrm{OH}_{\mathrm{k}_n})$$ onto a
cb-complemented subspace with constants independent of $n$ and
$$\big\langle \sigma(x), \omega(y) \big\rangle = \frac{1}{\ \mathrm{k}_n}
\sum_{j=1}^{\mathrm{k}_n} \mbox{tr}_{\mathcal{A}_n} (x_j^*y_j) =
\mbox{tr}_n(x^*y) = \langle x,y \rangle.$$ Therefore, the stated
properties of $\omega$ follow from those of the mapping $\sigma$.
\fin

Now we give a more explicit description of
$\mathcal{K}_{1,2}(\psi_n)$. Using the terminology introduced
before Lemma \ref{Lemma-PredualMAMS11}, the dual of the space
$\mathcal{L}_\infty^n$ is given by the following direct sum
$$\mathcal{L}_1^n = \big( R_n \otimes_h C_n \big)
\oplus_2 \big( R_n \otimes_h \mathrm{OH}_n \big) \oplus_2 \big(
\mathrm{OH}_n \otimes_h C_n \big) \oplus_2 \big( \mathrm{OH}_n
\otimes_h \mathrm{OH}_n \big).$$ Thus, we may consider the map
$$\Psi_n: \mathcal{L}_1^n \to L_1(q_n \mathcal{B}(\ell_2) q_n)$$
given by
$$\Psi_n (x_1, x_2, x_3,x_4) = x_1 + x_2
d_{\psi_n}^{\frac{1}{4}} + d_{\psi_n}^{\frac{1}{4}} x_3 +
d_{\psi_n}^{\frac{1}{4}} x_4 d_{\psi_n}^{\frac{1}{4}}.$$ Then it
is easily checked that $\ker \Psi_n =
\mathcal{J}_{\infty,2}(\psi_n)^{\perp}$ with respect to the
anti-linear duality bracket and we deduce
$\mathcal{K}_{1,2}(\psi_n) = \mathcal{L}_1^n / \ker \Psi_n$. The
finite-dimensional spaces defined so far allow us to take direct
limits
$$\mathcal{J}_{\infty,2} (\psi) = \overline{\bigcup_{n \ge 1}
\mathcal{J}_{\infty,2}(\psi_n)} \quad \mbox{and} \quad
\mathcal{K}_{1,2}(\psi) = \overline{\bigcup_{n \ge 1}
\mathcal{K}_{1,2}(\psi_n)}.$$

\begin{lemma} \label{Lemma-Direct-Sum11}
Let $\lambda_1, \lambda_2, \ldots \in \R_+$ be a sequence of
strictly positive numbers and define the diagonal operator
$\mathsf{d}_\lambda = \sum_k \lambda_k e_{kk}$ on $\ell_2$. Let us
equip the space $graph(\mathsf{d}_\lambda)$ with the following
operator space structures $$\begin{array}{rclcl} R \cap
\ell_2^{oh}(\lambda) & = & graph(\mathsf{d}_\lambda) & \subset & R
\oplus_2 \mathrm{OH}, \\ [3pt] C \cap \ell_2^{oh}(\lambda) & = &
graph(\mathsf{d}_\lambda) & \subset & C \hskip0.5pt \oplus_2
\mathrm{OH}. \end{array}$$ Then, if we consider the dual spaces
$$\begin{array}{rclcl} C + \ell_2^{oh}(\lambda) & = & \big( C
\oplus_2 \mathrm{OH} \big) \big/ \big( R \cap \ell_2^{oh}(\lambda)
\big)^\perp,
\\ [3pt] R + \ell_2^{oh}(\lambda) & = & \big( R \oplus_2
\mathrm{OH} \big) \big/ \big( C \cap \ell_2^{oh}(\lambda)
\big)^\perp, \end{array}$$ there exists a n.s.s.f. weight $\psi$
on $\mathcal{B}(\ell_2)$ such that
$$\big( R + \ell_2^{oh}(\lambda) \big) \otimes_h
\big( C + \ell_2^{oh}(\lambda) \big) = \mathcal{K}_{1,2}(\psi).$$
\end{lemma}

\dem Let $q_n$ be the projection $$\sum_{k=1}^\infty \alpha_k
\hskip1pt \delta_k \in \ell_2 \mapsto \sum_{k=1}^n \alpha_k
\hskip1pt \delta_k \in \ell_2^n.$$ If $\widehat{q}_n = q_n \oplus
q_n$, we define the subspaces
$$\begin{array}{lclcl} q_n \big( C + \ell_2^{oh}(\lambda) \big)
\!\!\!\! & = & \!\!\!\! \Big\{ \widehat{q}_n(a,b) + \big( R \cap
\ell_2^{oh}(\lambda) \big)^\perp \, \big| \ (a,b) \in C \oplus_2
\mathrm{OH} \Big\} \!\!\!\! & \subset & \!\!\!\! C +
\ell_2^{oh}(\lambda), \\ [5pt] q_n \big( R + \ell_2^{oh}(\lambda)
\big) \!\!\!\! & = & \!\!\!\! \Big\{ \widehat{q}_n(a,b) + \big( C
\cap \ell_2^{oh}(\lambda) \big)^\perp \, \big| \ (a,b) \in R
\oplus_2 \mathrm{OH} \Big\} \!\!\!\! & \subset & \!\!\!\! R +
\ell_2^{oh}(\lambda).
\end{array}$$ Note that, since the corresponding annihilators are
$q_n$-invariant, these are quotients of $C_n \oplus_2
\mathrm{OH}_n$ and $R_n \oplus_n \mathrm{OH}_n$ respectively.
Moreover, recalling that $q_n(x) \to x$ as $n \to \infty$ in the
norms of $R, \mathrm{OH}, C$, it is not difficult to see that we
may write the Haagerup tensor product $\big( R +
\ell_2^{oh}(\lambda) \big) \otimes_h \big( C +
\ell_2^{oh}(\lambda) \big)$ as the direct limit
$$\overline{\bigcup_{n \ge 1} q_n(R + \ell_2^{oh}(\lambda))
\otimes_h q_n(C + \ell_2^{oh}(\lambda))}.$$ Therefore, it suffices
to show that $$q_n(R + \ell_2^{oh}(\lambda)) \otimes_h q_n(C +
\ell_2^{oh}(\lambda)) = \mathcal{K}_{1,2}(\psi_n),$$ where
$\psi_n$ denotes the restriction to $q_n \mathcal{B}(\ell_2) q_n$
of some \emph{n.s.s.f.} weight $\psi$. However, by duality this is
equivalent to see that $q_n(C \cap \ell_2^{oh}(\lambda)) \otimes_h
q_n(R \cap \ell_2^{oh}(\lambda)) = \mathcal{J}_{\infty,2}(\psi_n)$
where the spaces $q_n(R \cap \ell_2^{oh}(\lambda))$ and $q_n(C
\cap \ell_2^{oh}(\lambda))$ are the span of $$\Big\{ (\delta_k,
\lambda_k \delta_k) \, \big| \ 1 \le k \le n \Big\}$$ in $R_n
\oplus_2 \mathrm{OH}_n$ and $C_n \oplus_2 \mathrm{OH}_n$
respectively. Indeed, we have
\begin{eqnarray*}
q_n \big( C + \ell_2^{oh}(\lambda) \big) & = & \big( C_n \oplus_2
\mathrm{OH}_n \big) / q_n(R \cap \ell_2^{oh}(\lambda))^{\perp},
\\ q_n \big( R + \ell_2^{oh}(\lambda) \big) & = & \big(
R_n \oplus_2 \mathrm{OH}_n \big) / q_n(C \cap
\ell_2^{oh}(\lambda))^{\perp},
\end{eqnarray*}
completely isometrically. Using row/column terminology in terms of
matrix units
\begin{eqnarray*}
q_n \big( C \cap \ell_2^{oh}(\lambda) \big) & = & \mbox{span}
\Big\{ (e_{i1} \hskip0.5pt , \hskip0.5pt \lambda_i \hskip0.5pt
e_{i1} \hskip0.5pt ) \in C_n \oplus_2
\hskip1pt \mathrm{OH}_n \Big\}, \\
q_n \big( R \cap \ell_2^{oh}(\lambda) \big) & = & \mbox{span}
\Big\{ (e_{1j}, \lambda_j e_{1j}) \in R_n \oplus_2 \mathrm{OH}_n
\Big\}.
\end{eqnarray*}
Therefore, the space $q_n(C \cap \ell_2^{oh}(\lambda)) \otimes_h
q_n(R \cap \ell_2^{oh}(\lambda))$ is the subspace
$$\mbox{span} \Big\{ (e_{ij}, \lambda_j e_{ij}, \lambda_i e_{ij},
\lambda_i \lambda_j e_{ij} ) \Big\} = \Big\{ (z, z
\mathsf{d}_{\lambda}, \mathsf{d}_\lambda z, \mathsf{d}_\lambda z
\mathsf{d}_\lambda) \, \big| \ z \in q_n \mathcal{B}(\ell_2) q_n
\Big\}$$ of the space $\mathcal{L}_\infty^n$ defined above. Then,
we define $\gamma_k \in \R_+$ by the relation $\lambda_k =
\gamma_k^{\frac14}$ and consider the \emph{n.s.s.f.} weight $\psi$
on $\mathcal{B}(\ell_2)$ induced by $\mathsf{d}_\gamma$. In
particular, we immediately obtain
$$q_n \big( C \cap \ell_2^{oh}(\lambda) \big) \otimes_h
q_n \big( R \cap \ell_2^{oh}(\lambda) \big) = \Big\{ \big( z, z
d_{\psi_n}^{\frac{1}{4}}, d_{\psi_n}^{\frac{1}{4}} z,
d_{\psi_n}^{\frac{1}{4}} z d_{\psi_n}^{\frac{1}{4}} \big)
\Big\}.$$ The space on the right is by definition
$\mathcal{J}_{\infty,2}(\psi_n)$. This completes the proof. \fin

\subsection{Proof of Theorem C}

In this last paragraph, we prove Theorem C. Let us begin with some
preliminary results. The next lemma has been known to Xu and the
first-named author for quite some time. We refer to Xu's paper
\cite{X3} for an even more general statement than the result
presented below.

\begin{lemma} \label{Lemma-Xu-Partition}
Given $1 \le p \le \infty$ and a closed subspace $\mathrm{X}$ of
$R_p \oplus_2 \mathrm{OH}$, there exist closed subspaces
$\mathcal{H}_1, \mathcal{H}_2, \mathcal{K}_1, \mathcal{K}_2$ of
$\ell_2$ and an injective closed densely-defined operator
$\Lambda: \mathcal{K}_1 \to \mathcal{K}_2$ with dense range such
that $$\mathrm{X} \simeq_{cb} \mathcal{H}_{1,r_p} \oplus_2
\mathcal{H}_{2,oh} \oplus_2 graph(\Lambda),$$ where the graph of
$\Lambda$ is regarded as a subspace of $\mathcal{K}_{1,r_p}
\oplus_2 \mathcal{K}_{2,oh}$ and the relevant constants in the
complete isomorphism above do not depend on the subspace
$\mathrm{X}$. Moreover, since $R_p = C_{p'}$ the same result can
be written in terms of column spaces.
\end{lemma}

In the following, we shall also need to recognize Pisier's
operator Hilbert space $\mathrm{OH}$ as the graph of certain
diagonal operator on $\ell_2$. More precisely, the following
result will be used below.

\begin{lemma} \label{Lemma-OH-graph}
Given $1 \le p \le \infty$, there exists a sequence $\lambda_1,
\lambda_2, \ldots$ in $\mathbb{R}_+$ for which the associated
diagonal map $\mathsf{d}_\lambda = \sum_k \lambda_k e_{kk}: R_p
\to \mathrm{OH}$ satisfies the following complete isomorphism
$$\mathrm{OH} \simeq_{cb} graph (\mathsf{d}_\lambda).$$
\end{lemma}

\dem Let us define
$$u: \delta_k \in \mathrm{OH} \mapsto (\lambda_k^{-1} \delta_k,
\delta_k) \in graph(\mathsf{d}_\lambda).$$ The mapping $u$
establishes a linear isomorphism between $\mathrm{OH}$ and
$graph(\mathsf{d}_\lambda)$. The inverse map of $u$ is the
coordinate projection into the second component, which is clearly
a complete contraction. Regarding the cb-norm of $u$, since
$graph(\mathsf{d}_\lambda)$ is equipped with the o.s.s. of $R_p
\oplus_2 \mathrm{OH}$, we have
$$\|u\|_{cb} = \sqrt{1 + \xi^2}$$ with $\xi$ standing for the
cb-norm of $\mathsf{d}_{\lambda^{-1}}: \mathrm{OH} \to R_p$. We
claim that $$\xi \le \Big( \summ_k |\lambda_k^{-1}|^4
\Big)^{\frac14},$$ so that it suffices to take $\lambda_1,
\lambda_2, \ldots$ large enough to deduce the assertion. Indeed,
it is well-known that the inequality above holds for the map
$\mathsf{d}_{\lambda^{-1}}: \mathrm{OH} \to R$ and also for
$\mathsf{d}_{\lambda^{-1}}: \mathrm{OH} \to C$. Therefore, our
claim follows by complex interpolation. \fin

\begin{remark}
\emph{The constants in Lemma \ref{Lemma-OH-graph} are uniformly
bounded on $p$.}
\end{remark}

Before beginning with the proof, we need a bit more preparation.
The following discretization result might be also well-known.
Nevertheless, since we are not aware of any reference for it, we
include the proof for the sake of completeness.

\begin{lemma} \label{Lemma-Diagonalization}
Given $1 \le p \le \infty$ and a closed densely-defined operator
$\Lambda: R_p \to \mathrm{OH}$ with dense range in $\mathrm{OH}$,
there exists a diagonal operator $\mathsf{d}_{\lambda} = \sum_k
\lambda_k e_{kk}$ on $\ell_2$ such that, when regarded as a map
$\mathsf{d}_\lambda: R_p \to \mathrm{OH}$, we obtain
$$graph(\mathsf{d}_\lambda) \simeq_{cb} graph(\Lambda).$$
Moreover, the relevant constants in the cb-isomorphism above do
not depend on $\Lambda$.
\end{lemma}

\dem Let us first assume that $\Lambda$ is positive. Then, since
$R_p$ is separable we deduce from spectral calculus \cite{KR} that
there exists a $\sigma$-finite measure space $(\Omega,
\mathcal{F}, \mu)$ for which $\Lambda$ is unitarily equivalent to
a multiplication operator $M_f: L_2(\Omega) \to L_2(\Omega)$. Thus
we may assume $\Lambda=M_f$. Now, we employ a standard procedure
to create a diagonal operator. Given $\delta>0$, we may
approximate the function $f$ by $$g=\summ_k (k\delta)
1_{k\delta<f\le (k+1)\delta}.$$ This yields a $1+\delta$
cb-isomorphism $$graph (M_f) \simeq_{cb} graph(M_g).$$ Therefore,
defining the measurable sets
$$\Omega_k = \Big\{ w \in \Omega \, \big| \ k \delta< f(w) \le
(k+1)\delta \Big\},$$ we have that $L_2(\Omega_k)$ is isomorphic
to $\ell_2(n_k)$ with $0 \le n_k = \dim L_2(\Omega_k) \le \infty$.
Choosing an orthonormal basis for $L_2(\Omega_k)$, we find that
$M_g$ is similar to $\mathsf{d}_{\lambda}$ where $\lambda_k = k
\delta$ with multiplicity $n_k$. This gives the assertion for
positive operators. If $\Lambda$ is not positive, we consider the
polar decomposition $\Lambda = u |\Lambda|$. By extension we may
assume that $u$ is a unitary. Thus, we get a cb-isometry $graph
(\Lambda) \simeq_{cb} graph (|\Lambda|)$ and the general case can
be reduced to the case of positive operators. \fin

\begin{lemma}\label{VNA-OH} Let $\mathcal{M}$ be a von Neumann
algebra. Then \[ L_1(\mathcal{M}; \mathrm{OH}) = [L_1(\mathcal{M};
C),L_1(\mathcal{M},R)]_{\frac12} \] completely embeds into
$L_1(\mathcal{A})$ for some von Neumann algebra $\mathcal{A}$.
\hskip-3pt Moreover, we have
\begin{itemize}
\item[\textbf{i)}] If $\mathcal{M}$ is $\mathrm{QWEP}$, we can
choose $\mathcal{A}$ to be $\mathrm{QWEP}$.

\item[\textbf{ii)}] If $\mathcal{M}$ is hyperfinite, we can choose
$\mathcal{A}$ to be hyperfinite.
\end{itemize}
\end{lemma}

\dem We recall from Pisier's theorem \cite{P0} (see also
\cite{JP2,X}) that
\[ \Big\| \summ_k x_k \otimes \delta_k \Big\|_{L_1(\mathcal{M};
\mathrm{OH})} = \inf_{x_k = a y_k b} \|a\|_4 \, \big( \summ_k
\|y_k\|_2^2)^{\frac12} \, \|b\|_4 . \] The first part of the proof
recaptures Pisier's argument in \cite{P4}. Pisier's exercise with
endpoints $L_1(\mathcal{M};C)$ and $L_1(\mathcal{M};R)$ implies
that $L_1(\mathcal{M}; \mathrm{OH})$ is completely isomorphic to
the quotient of ${\mathcal F}(L_1(\mathcal{M};C), L_1(\mathcal{M};
R))$ by the kernel of the mapping $$\mathcal{Q}: f \in {\mathcal
F}(L_1(\mathcal{M};C), L_1(\mathcal{M}; R)) \mapsto f(1/2) \in
L_1(\mathcal{M}; \mathrm{OH}).$$ Here ${\mathcal
F}(L_1(\mathcal{M};C), L_1(\mathcal{M}; R))$ is viewed as a
subspace of $$L_1(\mathcal{M}; L_2^c(\partial_0,\ell_2)) \oplus_1
L_1(\mathcal{M}; L_2^r(\partial_0,\ell_2))$$ and we find a
completely isomorphic embedding \[ L_1(\mathcal{M}; \mathrm{OH})
\subset L_1(\mathcal{M}; L_2^c(\partial_0,\ell_2)) \oplus_1
L_1(\mathcal{M}; L_2^r(\partial_1,\ell_2)) / \mathrm{ker}
\mathcal{Q} . \] It is worth mentioning that formally we might
need a finite $\mathcal{M}$ here to make the interpolation
argument work. However, this is no restriction in view of
Haagerup's reduction procedure \cite{H2}. Write $\mathrm{ker}
\mathcal{Q}$ as a tensor product $L_1(\mathcal{M}) \otimes
\mathrm{ker} \, q$ with $$q: f \in \mathcal{F}(C,R) \mapsto f(1/2)
\in \mathrm{OH}$$ and diagonalize $\mathrm{ker} \, q$ as in Lemma
\ref{Lemma-Xu-Partition}. Hence, it suffices to consider the
quotient space
\begin{eqnarray*}
\lefteqn{\hskip-20pt \Big\| \summ_k x_k \otimes \delta_k
\Big\|_{L_1(\mathcal{M}; Q(\la,\mu))}} \\ & = & \inf_{x_k = a_{kj}
+ b_{kj}} \Big\| \big( \summ_{k,j} \la_j \, a_{kj}^* a_{kj}
\big)^{\frac12} \Big\|_1 + \Big\| \big( \summ_{k,j} \mu_j \,
b_{kj} b_{kj}^* \big)^{\frac12} \Big\|_1 .
\end{eqnarray*}
Indeed, a continuous version of this formula has been obtained in
\cite{J2} and the proof there generalizes to arbitrary von Neumann
algebras. Now, we may apply Pisier's approach \cite{P4} and find
an embedding of $L_1(\mathcal{M}; Q(\la,\mu))$ in $L_1(\mathcal{M}
\ten \Gamma(\la,\mu))$ where $\Gamma(\la,\mu)$ is the free
quasi-free factor introduced by Shlyakhtenko.

\vskip5pt

To preserve hyperfiniteness, something that will not be applied in
this paper but in \cite{JP4}, we have to refer to \cite{J3}. The
results there are also stated in the operator space setting.
However, the isomorphism is based on the central limit procedure
and a Khintchine type inequality which holds in full generality.
For the central limit procedure we shall first consider finitely
many coordinates $x_1, x_2, \ldots, x_m$ and finite sequences
$$\la(n) = (\la_1, \la_2, \ldots, \la_n) \quad \mbox{and} \quad
\mu(n) = (\mu_1, \mu_2, \ldots, \mu_n).$$  Let $R(\la, \mu)$
denote the Araki-Woods factor constructed from $(\la,\mu)$ with
quasi-free state $\phi$ and density $d$. Let $R(\la(n),\mu(n))$ be
the corresponding finite dimensional matrix algebras and
$$\mathcal{R}_{m,n} = \bigotimes_{1 \le k \le m}
R(\la(n),\mu(n)).$$ This yields elements $\xi_k \in
L_1(\mathcal{R}_{m,n})$ such that
\[ \Big\| \summ_{k \le m} x_k \otimes \delta_k
\Big\|_{L_1(\mathcal{M}; Q(\la(n),\mu(n)))} \sim_c \Big\| \summ_{k
\le m} x_k \ten \xi_k(\la_n,\mu_n) \Big\|_{L_1(\mathcal{M}
\bar\otimes \mathcal{R}_{m,n})} .\] Then we take first $n\to
\infty$ and then $m\to \infty$. Here we use the isometric
embedding of $L_1(\mathcal{M}, Q(\la,\mu))$ in its bidual and the
fact that finite sequences are dense in $L_1(\mathcal{M};
\mathrm{OH})$. We obtain an infinite tensor product $\mathcal{R} =
\ten_{k \ge 1} R(\la,\mu)$. This yields an embedding in
$L_1(\mathcal{M} \bar\ten \mathcal{R})$ where $\mathcal{R}$ is an
Araki-Woods factor. Indeed, the $\mathrm{III}_1$ factor will do.
By replacing $\mathcal{M}$ with $M_m \ten \mathcal{M}$ we see that
this embedding is automatically a complete isomorphism. Thus for
$\mathcal{M}$ hyperfinite (resp. $\mathrm{QWEP}$) the tensor
product $\mathcal{M} \bar\ten \mathcal{R}$ has the same property.
\fin

\demC By injectivity of the Haagerup tensor product, we may assume
that $(\mathrm{X}_1, \mathrm{X}_2) \in \mathcal{Q}(R \oplus_2
\mathrm{OH}) \times \mathcal{Q}(C \oplus_2 \mathrm{OH})$. In
particular, the duals $\mathrm{X}_1^*$ and $\mathrm{X}_2^*$ are
subspaces of $C \oplus_2 \mathrm{OH}$ and $R \oplus_2 \mathrm{OH}$
respectively. Therefore, according to Lemma
\ref{Lemma-Xu-Partition} above, we may find Hilbert spaces
$\mathcal{H}_{ij}$ and $\mathcal{K}_{ij}$ for $i,j=1,2$ such that
\begin{eqnarray*}
\mathrm{X}_1^* & \simeq_{cb} & \mathcal{H}_{11,c} \oplus_2
\mathcal{H}_{12,oh} \oplus_2 graph(\Lambda_1), \\
\mathrm{X}_2^* & \simeq_{cb} & \mathcal{H}_{21,r} \oplus_2
\mathcal{H}_{22,oh} \oplus_2 graph(\Lambda_2),
\end{eqnarray*}
where the operators $\Lambda_1: \mathcal{K}_{11,c} \to
\mathcal{K}_{12,oh}$ and $\Lambda_2: \mathcal{K}_{21,r} \to
\mathcal{K}_{22,oh}$ are injective, closed, densely-defined with
dense range. On the other hand, using the complete isometries
$\mathcal{H}_r^* = \mathcal{H}_c$ and $\mathcal{H}_c^* =
\mathcal{H}_r$, we easily obtain the cb-isomorphisms
\begin{eqnarray*}
\mathrm{X}_1 & \simeq_{cb} & \mathcal{H}_{11,r} \oplus_2
\mathcal{H}_{12,oh} \oplus_2 \Big( \big( \mathcal{K}_{11,r}
\oplus_2 \mathcal{K}_{12,oh} \big) \big/ graph(\Lambda_1)^\perp
\Big), \\ \mathrm{X}_2 & \simeq_{cb} & \mathcal{H}_{21,c} \oplus_2
\mathcal{H}_{22,oh} \oplus_2 \Big( \big( \mathcal{K}_{21,c}
\oplus_2 \mathcal{K}_{22,oh} \big) \big/ graph(\Lambda_2)^\perp
\Big).
\end{eqnarray*}
Let us set for the sequel
\begin{eqnarray*}
\mathcal{Z}_1 & = & \big( \mathcal{K}_{11,r} \oplus_2
\mathcal{K}_{12,oh} \big) \big/ graph(\Lambda_1)^\perp, \\
\mathcal{Z}_2 & = & \big( \mathcal{K}_{21,c} \oplus_2
\mathcal{K}_{22,oh} \big) \big/ graph(\Lambda_2)^\perp.
\end{eqnarray*}
Then, we have the following cb-isometric inclusion
\begin{eqnarray} \label{Eq-6terms11}
\mathrm{X}_1 \otimes_h \mathrm{X}_2 \subset \mathrm{A}_1 \oplus_2
\mathrm{A}_2 \oplus_2 \mathrm{A}_3 \oplus_2 \mathrm{A}_4 \oplus_2
\mathrm{A}_5 \oplus_2 \mathrm{A}_6,
\end{eqnarray}
where the $\mathrm{A}_j$'s are given by
\begin{eqnarray*}
\mathrm{A}_1 & = & \mathcal{Z}_1 \otimes_h \mathcal{Z}_2 \\
\mathrm{A}_2 & = & \mathcal{H}_{11,r} \otimes_h \mathrm{X}_2
\\ \mathrm{A}_3 & = & \mathrm{X}_1 \otimes_h \mathcal{H}_{21,c} \\
\mathrm{A}_4 & = & \mathcal{H}_{12,oh} \otimes_h \mathcal{Z}_2
\\ \mathrm{A}_5 & = & \mathcal{Z}_1 \otimes_h
\mathcal{H}_{22,oh} \\ \mathrm{A}_6 & = & \mathcal{H}_{12,oh}
\otimes_h \mathcal{H}_{22,oh}.
\end{eqnarray*}
We now reduce the proof to the construction of a cb-embedding
$\mathcal{Z}_1 \otimes_h \mathcal{Z}_2 \to L_1(\mathcal{A})$ for
some $\mathrm{QWEP}$ von Neumann algebra $\mathcal{A}$. Indeed,
according to \cite{J2} we know that $\mathrm{OH}$ cb-embeds in
$L_1(\mathcal{A})$ for some $\mathrm{QWEP}$ type $\mathrm{III}$
factor $\mathcal{A}$. Hence, the term $\mathrm{A}_6$ automatically
satisfies the assertion. A similar argument works for the terms
$\mathrm{A}_2$ and $\mathrm{A}_3$. Indeed, they clearly embed into
the vector-valued Schatten classes $S_1(\mathrm{X}_2)$ and
$S_1(\mathrm{X}_1)$ cb-isometrically. On the other hand, since
$\mathrm{OH} \in \mathcal{QS}(C \oplus R)$ by \lq\lq Pisier's
exercise\rq\rq${}$ and we have by hypothesis
$$\mathrm{X}_1 \in \mathcal{QS}(R \oplus_2 \mathrm{OH}) \quad
\mbox{and} \quad \mathrm{X}_2 \in \mathcal{QS}(C \oplus_2
\mathrm{OH}),$$ both $\mathrm{X}_1$ and $\mathrm{X}_2$ are
cb-isomorphic to an element in $\mathcal{QS}(C \oplus R)$.
According to \cite{J2} one more time, we know that any operator
space in $\mathcal{QS}(C \oplus R)$ cb-embeds into
$L_1(\mathcal{A})$ for some $\mathrm{QWEP}$ von Neumann algebra
$\mathcal{A}$. Thus, the spaces $S_1(\mathrm{X}_1)$ and
$S_1(\mathrm{X}_2)$ also satisfy the assertion. Finally, for
$\mathrm{A}_4$ and $\mathrm{A}_5$ we apply Lemma
\ref{Lemma-OH-graph} and write $\mathrm{OH}$ as the graph of a
diagonal operator on $\ell_2$. By the self-duality of
$\mathrm{OH}$ we conclude that these terms can be regarded as
particular cases of $\mathrm{A}_1$.

\vskip5pt

It remains to see that the term $\mathcal{Z}_1 \otimes_h
\mathcal{Z}_2$ satisfies the assertion. According to the
discretization Lemma \ref{Lemma-Diagonalization}, we may assume
that the graphs appearing in the terms $\mathcal{Z}_1$ and
$\mathcal{Z}_2$ above are graphs of diagonal operators
$\mathsf{d}_{\lambda_1}$ and $\mathsf{d}_{\lambda_2}$. In fact, by
polar decomposition as in the proof of Lemma
\ref{Lemma-Diagonalization}, we may also assume that both diagonal
operators are positive. Moreover, by adding a perturbation term we
can take the eigenvalues $\lambda_{1k}, \lambda_{2k} \in
\mathbb{R}_+$ strictly positive. Indeed, if we replace
$\lambda_{jk}$ by $\xi_{jk} = \lambda_{jk} + \varepsilon_k$ for
$j=1,2$, the new diagonal operators $\mathsf{d}_{\xi_1}$ and
$\mathsf{d}_{\xi_2}$ satisfy the cb-isomorphisms
$$graph(\mathsf{d}_{\lambda_j}) \simeq_{cb}
graph(\mathsf{d}_{\xi_j}) \quad \mbox{for} \quad j=1,2$$ where
(arguing as in Lemma \ref{Lemma-OH-graph} above) the cb-norms are
controlled by
$$\Big( \summ_k |\varepsilon_k|^4 \Big)^{\frac14}.$$
Therefore, taking the $\varepsilon_k$'s small enough, we may write
\begin{eqnarray*}
\mathcal{Z}_1 = \big( R \oplus_2 \mathrm{OH} \big) \big/ \big( C
\cap \ell_2^{oh}(\lambda_1) \big)^{\perp} = R +
\ell_2^{oh}(\lambda_1) & \mbox{with} &
\mathsf{d}_{\lambda_1}: C \to \mathrm{OH}, \\
\mathcal{Z}_2 = \big( C \oplus_2 \mathrm{OH} \big) \big/ \big( R
\cap \ell_2^{oh}(\lambda_2) \big)^{\perp} = C +
\ell_2^{oh}(\lambda_2) & \mbox{with} & \mathsf{d}_{\lambda_2}: R
\to \mathrm{OH},
\end{eqnarray*}
where $\mathsf{d}_{\lambda_1}, \mathsf{d}_{\lambda_2}$ are
positive and invertible. Now we set
$$\lambda_k = \lambda_{[k+1/2]}^j \quad  \mbox{for} \quad k \equiv
j \ (\mathrm{mod} \, 2).$$ If we define $\mathsf{d}_{\lambda} =
\sum_k \lambda_k e_{kk}$, we find a complete embedding
$$\mathcal{Z}_1 \otimes_h \mathcal{Z}_2 \subset \big( R +
\ell_2^{oh}(\lambda) \big) \ten_h \big( C + \ell_2^{oh}(\lambda)
\big).$$ According to Lemma \ref{Lemma-Direct-Sum11}, we can
regard $\mathcal{Z}_1 \otimes_h \mathcal{Z}_2$ as a subspace of
$$\mathcal{K}_{1,2}(\psi) = \overline{\bigcup_{n \ge 1}
\mathcal{K}_{1,2}(\psi_n)}^{\null},$$ for some \emph{n.s.s.f.}
weight $\psi$ on $\mathcal{B}(\ell_2)$. It remains to construct a
completely isomorphic embedding from $\mathcal{K}_{1,2}(\psi)$
into $L_1(\mathcal{A})$ for some $\mathrm{QWEP}$ algebra
$\mathcal{A}$. To that aim we assume that the numbers
$\mathrm{k}_n = \psi_n(q_n)$ are non-decreasing positive integers.
This can always be achieved by the same perturbation argument used
above. This will allow us to apply Lemma
\ref{Lemma-PredualMAMS11}. Now, in order to cb-embed
$\mathcal{K}_{1,2}(\psi)$ into $L_1(\mathcal{A})$ it suffices to
construct a cb-embedding of $\mathcal{K}_{1,2}(\psi_n)$ into
$L_1(\mathcal{A}_n')$ for some $\mathcal{A}_n'$ being
$\mathrm{QWEP}$ and with relevant constants independent of $n$.
Indeed, if so we may consider an ultrafilter $\mathcal{U}$
containing all the intervals $(n,\infty)$, so that we have a
completely isometric embedding
$$\mathcal{K}_{1,2}(\psi) = \overline{\bigcup_{n \ge 1}
\mathcal{K}_{1,2}(\psi_n)}^{\null} \to \prodd_{n, \mathcal{U}}
\mathcal{K}_{1,2}(\psi_n).$$ Then, according to \cite{Ra}, our
assumption provides a cb-embedding
$$\mathcal{K}_{1,2}(\psi) \to L_1(\mathcal{A}) \quad
\mbox{with} \quad \mathcal{A} = \Big( \prodd_{n, \mathcal{U}}
{\mathcal{A}_n'}_* \Big)^*.$$ Moreover, we know from \cite{J5}
that $\mathcal{A}$ is $\mathrm{QWEP}$ provided the
$\mathcal{A}'_n$'s are. Thus, it remains to construct the
cb-embeddings $\mathcal{K}_{1,2}(\psi_n) \to L_1(\mathcal{A}'_n)$.
According to Lemma \ref{Lemma-PredualMAMS11}, we have a
cb-embedding
$$\mathcal{K}_{1,2}(\psi_n) \to L_1(\mathcal{A}_n;
\mathrm{OH}_{\mathrm{k}_n}).$$ Now we observe that $\mathcal{A}_n$
is QWEP since it is the free product of $\mathrm{k}_n$ copies of
$\mathrm{M}_n \oplus \mathrm{M}_n$ and we know from \cite{J2} that
the $\mathrm{QWEP}$ is stable under reduced free products. Hence
Lemma \ref{VNA-OH} implies that $L_1(\mathcal{A}_n;
\mathrm{OH}_{\mathrm{k}_n})\subset L_1(\mathcal{A}'_n)$ such that
$\mathcal{A}'_n$ is also QWEP. \fin

\begin{corollary} \label{SqQWEP11}
If $1 < q \le 2$, we have $$S_q \hookrightarrow_{cb}
L_1(\mathcal{A})$$ for some von Neumann algebra $\mathcal{A}$
satisfying the $\mathrm{QWEP}$.
\end{corollary}

\dem Using the complete isometry $$S_q = C_q \otimes_h R_q,$$ the
assertion follows combining Lemma \ref{Lemma-Motivation} and
Theorem \ref{QSQSQS}. \fin

\bibliographystyle{amsplain}

\vskip20pt

\noindent \textbf{Marius Junge} \\
\textsc{Department of Mathematics} \\
\textsc{University of Illinois at Urbana-Champaign} \\ 273 Altgeld
Hall, 1409 W. Green Street, Urbana, IL 61801, USA \\
\texttt{junge@math.uiuc.edu}

\vskip0.3cm

\noindent \textbf{Javier Parcet} \\
\textsc{Departamento de Matem{\'a}ticas} \\
\textsc{Instituto de Matem}{\scriptsize {\'A}}\textsc{ticas
y F}{\scriptsize {\'I}}\textsc{sica Fundamental} \\
\textsc{Consejo Superior de Investigaciones Cient}{\scriptsize
{\'I}}\textsc{ficas} \\ Depto de
Matem{\'a}ticas, Univ. Aut{\'o}noma de Madrid, 28049, Spain \\
\texttt{javier.parcet@uam.es}

\end{document}